\documentclass{amsart}
\usepackage{amssymb}
\usepackage{amsmath}
\usepackage{amsfonts}
\usepackage{latexsym}

\usepackage{color}
\usepackage[applemac]{inputenc}
\usepackage{enumerate}

 \newtheorem{definition}{Definition}[section]
 \newtheorem{remark}[definition]{Remark}
 
\theoremstyle{plain}      
 \newtheorem{proposition}[definition]{Proposition}
 \newtheorem{theorem}[definition]{Theorem}
 \newtheorem{corollary}[definition]{Corollary}
 \newtheorem{lemma}[definition]{Lemma}

\newcommand{\Tr}{\operatorname{Tr}}

\newcommand{\Hom}{\operatorname{Hom}}
\newcommand{\GL}{\operatorname{GL}}

\newcommand{\PGL}{\operatorname{PGL}}
\newcommand{\Sym}{\operatorname{Sym}}
\newcommand{\Conf}{\operatorname{Conf}}
\newcommand{\Met}{\operatorname{Met}}
\newcommand{\PMod}{\operatorname{PMod}}
\newcommand{\Mod}{\operatorname{Mod}}
\newcommand{\Aut}{\operatorname{Aut}}
\newcommand{\Out}{\operatorname{Out}}
\newcommand{\POut}{\operatorname{POut}}
\newcommand{\Inn}{\operatorname{Inn}}
\newcommand{\SE}{\operatorname{SE}}

\newcommand{\Sim}{\operatorname{Sim}}
\newcommand{\Diff}{\operatorname{Diff}}
\newcommand{\PDiff}{\operatorname{PDiff}}

\newcommand{\R}{{\mathbb R}}
\newcommand{\C}{{\mathbb C}}

\renewcommand{\b}{\beta}

\newcommand{\tq}{\, \big| \, }
\newcommand{\Tau}{\mathcal{T}}

\makeindex

\begin{document}

\title{On the Moduli Space of Singular Euclidean Surfaces}

\author{Marc Troyanov}  

\address{
 \'Ecole Polytechnique F{\'e}derale de
Lausanne (EPFL), \\
Institut de G{\'e}om{\'e}trie, Alg{\`e}bre et Topologie
\\ CH-1015 Lausanne - Switzerland
\\ email:\,\tt{marc.troyanov@epfl.ch}
}

\maketitle

\begin{abstract}
The goal of this paper is to develop  some aspects of the deformation
theory of piecewise flat structures on surfaces and use this theory
to construct new geometric structures on the moduli space of Riemann
surfaces.
\end{abstract}

\tableofcontents

\section*{Introduction}

\  The Teichmüller space of a punctured surface is the space of hyperbolic metrics with cusps up to isotopy on that surface, however, it can also be seen as the space of flat metrics with conical singularities of prescribed angles at the punctures up to isotopy and rescaling. The aim of the present paper is to use this fact and show how  the theory of piecewise flat surfaces and their deformations leads to new geometric structures in Teichmüller theory.

In the first section, which is rather elementary, we describe the
geometry of piecewise flat surfaces. The second section describes
the topology of punctured surfaces and their diffeomorphism
groups. In the third section we discuss the representation space
of a finitely generated group $\pi$ into the group $\SE(2)$ of
rigid motions in the euclidean plane.
In the last two sections, we apply the previous  results  to construct
a new geometric structure on the Riemann  moduli space
$\mathcal{M}_{g,n}$ of a surface $\Sigma$ of genus $g$ with $n$
punctures. More specifically, we show that this moduli space is
a good orbifold\footnote{Recall that an \emph{orbifold} is a space which is locally the quotient of a manifold by a finite group. A \emph{good orbifold} is globally the quotient of a manifold by a  group acting properly and discontinuously (but in general not freely).}  which admits a family of geometric structures locally 
modeled on the homogeneous spaces \ $ \Xi = \mathbb{T}^{2g}
\times  \mathbb{CP}^{2g+n-3}$.

\medskip

We now discuss our main result. We first define  a
\emph{punctured surface}\index{punctured surface}  $\Sigma_{g,n}$
of type $(g,n)$ to be a fixed connected closed orientable surface
$S$ of genus $g$ together with a distinguished set of $n$
pairwise distinct points $p_1,p_2,...,p_n \in \Sigma_{g,n}$. The
\emph{Teichm\"uller space}  $\mathcal{T}_{g,n}$ of $\Sigma_{g,n}$
is the set of  conformal structures on $\Sigma_{g,n}$ modulo
isotopies fixing the punctures (see section \ref{sec.teich} for a
precise definition). This space is a real analytic variety in a
natural way; if $2g-2 +n >0$, then it is isomorphic to $\R^{6g - 6
+ 2n}$. The group of orientation-preserving isotopy classes of diffeomorphisms of
$\Sigma_{g,n}$ fixing the punctures is called the \emph{pure
mapping class group} and denoted by $ \PMod_{g,n}$. It acts in a
natural way on the Teichm\"uller space  $\mathcal{T}_{g,n}$.

\bigskip

We are now in a position to state  the main result:

\medskip

\textbf{Theorem} \
\emph{Given a punctured surface $\Sigma_{g,n}$
such that $2g+n-2> 0$,  we can construct a group homomorphism
$$
 \Phi : \PMod_{g,n} \to \mathcal{G} = \Aut (\mathbb{T}^{2g})\times \PGL_{2g+n-2}(\mathbb{C})
$$
and a $\Phi$-equivariant  local homeomorphism
$$
\mathcal{H} : \mathcal{T}_{g,n} \to  \Xi = \mathbb{T}^{2g} \times
\mathbb{CP}^{2g+n-3}.
$$}
 
\medskip

To say that $\mathcal{H}$ is $\Phi$-equivariant means that
$
  \mathcal{H} (A \mu) = \Phi(A) \cdot \mathcal{H}(\mu)
$
for all $A\in \PMod_{g,n}$ and $\mu \in \mathcal{T}_{g,n}$.

\medskip

The pair $(\mathcal{H},\Phi)$ depends on $n$ parameters
$\b_1,\b_2,\dots ,\b_n \in (-1,\infty )$ such that $\sum_{j=1}^n
\b_j =  2g-2$ and no $\b_i\in \mathbb{Z}$.

\medskip 

The \emph{moduli space}  $\mathcal{M}_{g,n}$ of $\Sigma_{g,n}$ is
the set of  conformal structures on $\Sigma_{g,n}$ modulo
diffeomorphisms fixing the punctures. It is the quotient of the
Teichm\"uller space by the pure mapping class group of
$\Sigma_{g,n}$; in other words $\mathcal{M}_{g,n}$  is a good
orbifold whose universal cover is  $\mathcal{T}_{g,n}$ and
fundamental group is  $\PMod_{g,n}$. In the geometric language
of $(G,X)$-strucures on manifolds and orbifolds (see \cite{CEG,choi2004,goldman88,thurston97}), this theorem says that \emph{we have constructed a family  of geometric structures on the orbifold $\mathcal{M}_{g,n}$ which is
modeled on the homogeneous space $\Xi = \mathbb{T}^{2g}
\times\mathbb{CP}^{2g+n-3}$. This family is parametrized by the
$\beta_i's$.}

\medskip

The composition of the map $\mathcal{H}$ in the  Theorem with the projection on 
the torus $\mathbb{T}^{2g}$ gives us a map
$\rho : \mathcal{T}_{g,n} \to  \mathbb{T}^{2g}$ called the \emph{character map}. It was proved by  W.A. Veech that this map is a real analytic submersion. Its fibers describe a foliation whose leaves carry a geometric structure locally modelled on the complex
projective space $\mathbb{CP}^{2g+n-3}$, see \cite{veech93} for proofs of these facts and a discussion of other related geometric structures on $\mathcal{T}_{g,n}$.

\medskip

The proof of this theorem is based on the following strategy: we
first show that the Teichm\"uller space can be seen as a
deformation space of  \emph{flat metrics} on  $\Sigma_{g,n}$
having \emph{conical singularities} of prescribed angles at the
punctures. We associate to such a metric a homomorphism, called
the \emph{holonomy} of the metric,  from the fundamental group of
the surface to the group $\SE(2)$ of direct isometries of the
euclidean plane. We then show that such a homomorphism can be seen
as a point in the variety $\Xi$. In brief, $\mathcal{H} :
\mathcal{T}_{g,n} \to  \Xi$ maps the isotopy class of a singular
flat metric to its holonomy representation.

\medskip

In the special case of the punctured sphere, a stronger form of this
theorem has been obtained by Deligne and Mostow   \cite{DM} using some
techniques from algebraic geometry and by Thurston  \cite{thurston98}
using an approach closer to ours.

\medskip

To conclude this introduction, let us stress that the importance of piecewise flat  metrics in Teichmüller theory is illustrated by the large number of papers dedicated to this subject. In addition to the work of Veech and Thurston already quoted, let us mention the contributions of  Rivin \cite{rivin}, Bowditch \cite{Bow}, Epstein and Penner \cite{EpsPenn} to name a few. Piecewise euclidean metrics also appear in  quantum gravity and in topological quantum field theory, see \cite{ACM,CDM} and the references therein. Although the present paper starts with elementary considerations, the reader ought not to consider it as a global survey of this vast subject.

\medskip

\textbf{Acknowledgments.}  I would like to thank Babak Modami and Fran\c{c}ois Fillastre for having carefully read the manuscript and for their comments. Finally, this paper is dedicated to the memory of Michel Matthey.


\section{Piecewise flat surfaces}

\subsection{Euclidean triangulation on a surface}

A \emph{piecewise flat surface} is a metric space obtained by properly gluing 
a stock of euclidean triangles in such a way that
whenever two triangles  meet
along an edge,   they are glued by an isometry along that edge. More precisely:
\begin{definition}
A  \emph{euclidean triangulation}\index{euclidean
triangulation}   of a surface $\Sigma$ is  a set of pairs $\Tau
= \{(T_{\alpha},f_{\alpha})\}_{\alpha \in A}$ where each
$T_{\alpha}$ is a compact subset of $\Sigma$ and $f_{\alpha}:
T_{\alpha} {\to} \R^2$ is a homeomorphism onto a non degenerate
triangle  $f_{\alpha}(T_{\alpha})$ in the euclidean plane $\R^2$.
A subset $e$ of $T_{\alpha}$ is an \emph{edge} if $f_{\alpha}(e)$
is an edge of $f_{\alpha}(T_{\alpha})$ and a point  $p$ of
$T_{\alpha}$ is a \emph{vertex} if its image under
$f_{\alpha}$ is a vertex of $f_{\alpha}(T_{\alpha})$.

\smallskip

The eucliden triangulation $\Tau$  is subject to
the following conditions:
\begin{enumerate}[i)]
    \item The triangles cover the surface: $\Sigma = \bigcup_{\alpha }T_{\alpha}$.
  \item If ${\alpha}\neq{\beta}$, then the intersection $T_{\alpha} \cap T_{\beta}$
   is either empty, or an edge or a vertex.
  \item If $T_{\alpha} \cap T_{\beta}\neq \emptyset$,
   then there is an element $g_{\alpha\beta} \in E(2)$ (= the group of isometries of the euclidean plane) 
   such that
  $f_{\alpha}=g_{\alpha\beta}f_{\beta}$ on that intersection.
  \end{enumerate}
\end{definition}

An element $(T_{\alpha},f_{\alpha})\in \Tau$ is called a
\emph{triangle} or a \emph{2-simplex} of the triangulations, we
often just denote it by $T_{\alpha}$. The vertices and edges  are
called 0- and 1-simplices respectively.

\medskip

Two euclidean triangulations $\Tau = \{(T_{\alpha},f_{\alpha})\}_{\alpha \in A}$ and
$\Tau' = \{(T_{\alpha},f'_{\alpha})\}_{\alpha \in A}$ of the same surface $\Sigma$
are considered to be equal if they have the same simplices and, for any $\alpha \in A$,
there is an isometry $g_{\alpha} \in E(2)$ such that $f'_{\alpha}
= g_{\alpha}f_{\alpha}$.

\begin{definition}
A \emph{piecewise flat surface}\index{piecewise flat surface}   $(\Sigma,\mathcal{T})$  is a surface together with a euclidean triangulation.
\end{definition}

A piecewise flat surface  $(\Sigma,\mathcal{T})$ comes with a number of additional structures.
In particular there is a well defined \emph{area measure} which coincides with the 2-dimensional
Lesbegue measure on each euclidean triangle $T$. We can also define the
length $\ell (c)$ of an arbitrary curve $c : [0,1] \to \Sigma$ by the following axioms:
\begin{enumerate}[(i)]
  \item if $c$ is contained in a triangle  $T$ of $\mathcal{T}$, then $\ell (c)$ is
  the euclidean length.
  \item $\ell $ is additive: if $c$ is the concatenation of two curves $c_1c_2$, then $\ell (c)
  = \ell(c_1)+\ell(c_2)$.
\end{enumerate}
The piecewise flat surface is thus a length space (see \cite{bbi}
for this notion). If the surface is connected, then it is also a
metric space for the distance given by the infimum of the lengths
of all curves joining two given points.

\medskip

There is one more structure, called the \emph{singularity order}
and which  is defined as  the angular excess at the vertices
counted in number of turns. It tells us how singular each vertex
is  compared to an ordinary point; the precise definition is the
following:
\begin{definition}
The vertex $p\in \Sigma$ is said to be a  \emph{conical point of
total angle $\theta$} if
$$
 \theta = \sum_{j=1}^k\varphi_{j}
$$
where $\varphi_{1},\dots , \varphi_{k}$ are the angles of all the
triangles in ${\mathcal{T}}$ which are incident to $p$. The
\emph{singularity order} $\beta(p)$ of a vertex $p$ is the angular
excess at $p$ measured in number of turns:
$$
 \beta (p) =  \frac{\theta }{2\pi}- 1.
$$
We extend the function $\beta$ to all points of $\Sigma$ by
setting $\beta (x) = 0$ if $x\in \Sigma$ is not a vertex. The
point $x$ is then termed \emph{singular} if $\beta(x) \neq 0$
(i.e. if $\theta \neq 2\pi$) and \emph{regular} otherwise.
\end{definition}

\begin{proposition}[Gauss-Bonnet Formula]
For any  euclidean triangulation on a compact surface without boundary  $\Sigma$, we have
\begin{equation}\label{cond.GB}
 \chi (\Sigma) + \sum_{x\in\Sigma} \b(x) = 0 \, ,
 \end{equation}
where $\chi(\Sigma)$ is the Euler characteristic of the surface.
\end{proposition}

The proof is a direct counting argument based on  the  definition
of the Euler characteristic and the fact that the three  internal
angles of a euclidean triangle add up to $\pi$, see \cite{Tr86}.

\subsection{The  universal branched cover of a piecewise flat surface}

If  $(\Sigma, \mathcal{T})$ is a piecewise flat surface, we denote by $\Sigma' = \Sigma \setminus \{p_1,...,p_n \}$
the open surface obtained by removing the singular vertices $p_1,...,p_n \in \Sigma$.
\begin{definition}
A path $c : [0,1] \to \Sigma$   is
\emph{admissible}\index{admissible path} if it has finitely many
intersections with the edges of the triangulation and if $c(s) \in
\Sigma'$ for any $0< s <1$. A homotopy $c_t$  is an
\emph{admissible homotopy}\index{admissible homotopy} if  $s
\mapsto  c_t(s) \in \Sigma$ is an admissible path for any $0\leq
t \leq 1$.
\end{definition}

Let us choose a fixed triangle $T_0 \in \mathcal{T}$ and call it \emph{home} (or the \emph{base triangle}). We also choose a base point  $x_0$ in the interior of $T_0$.
\begin{definition}
The  \emph{universal branched cover}\index{universal branched cover} of $(\Sigma, \mathcal{T})$ is the euclidean two-dimensional complex  $\widehat{\mathcal{T}}$ obtained as follows: a $k$-simplex $\widehat{\sigma}$ of  $\widehat{\mathcal{T}}$, where $k=0,1$ or $2$,  is a pair
$(\sigma, [c])$ where $\sigma$ is a $k$-simplex of $\mathcal{T}$ and $[c]$ is an admissible homotopy class of paths  joining $T_0$ to a point in $\sigma$.
\end{definition}

The universal branched cover $\widehat{\mathcal{T}}$ is a
simplicial complex (which is not locally finite) and there is an
obvious simplicial map $\widehat{\mathcal{T}} \to \mathcal{T}$
sending $(\sigma, [c])$ to $\sigma$.

\medskip

We denote by $\widehat{\Sigma}$ the geometric  realization of
$\widehat{\mathcal{T}}$. This is a triangulated topological space
and it comes with a continuous surjective map $P:
\widehat{\Sigma} \to \Sigma$ sending each simplex of
$\widehat{\mathcal{T}}$ homeomorphically onto the corresponding
simplex in $\mathcal{T}$. We turn  $\widehat{\Sigma}$ into a
metric space (in fact a length space) by requiring $P$ to be an
isometry on each simplex (concretely, we give to each simplex
$\widehat{\sigma} = (\sigma, [c])$ in $\widehat{\mathcal{T}}$ the
geometry of the euclidean simplex $\sigma$ in $\mathcal{T}$).

\medskip

Another way to understand $\widehat{\Sigma}$ is the following: let
$\widetilde{\Sigma}'$ be the universal cover of $\Sigma'$. It is naturally
a length space (in fact a flat Riemannian surface) and $\widehat{\Sigma}$
is its metric completion.

\subsection{The development of a piecewise flat surface}
\label{sec.dev}

\begin{definition} An edge of the piecewise flat surface  $(\Sigma, \mathcal{T})$  is said to be \emph{interior} if it is not contained in the boundary of $\Sigma$.
The \emph{hinge}\index{hinge}   of an interior edge $e$   is the unique pair of triangles $T_1,T_2 \in \mathcal{T}$ which are incident with $e$.
\end{definition}

Given an interior edge $e$ with hinge $(T_1,T_2 )$
and an isometry $f_1 : T_1 \to \mathbb{R}^2$, there exists a unique
isometry $f_2 : T_2 \to \mathbb{R}^2$ such that  $f_1(T_1)$ and $f_2(T_2)$ have disjoint interiors and $f_1(e) = f_2(e)$. By juxtaposing these maps, we obtain a map
$$f_{e} =f_1 \cup f_2 :T_{1} \cup T_{2} \to \R^2$$
which is an isometry of the hinge onto a  quadrilateral in the
euclidean plane. The map $f_{e}$ just described is called an {\em
unfolding} of the hinge. One also says that $f_2$ is the
\emph{continuation} of $f_1$ across the edge $e$.

\medskip

The notions of hinge, unfolding and continuation of an isometry
across an edge are  similarly defined  on the universal branch
cover $\widehat{\Sigma}$.

\medskip

\begin{proposition}
 Let $(\Sigma, \mathcal{T})$ be a piecewise flat surface with home triangle $T_0$ and choose
 an isometry $f_0$ from  $T_0$ onto a triangle in  $\R^2$. Then there exists a unique map
 $f : \widehat{\Sigma} \to \R^2$ such that $f$ coincides with $f_0$ on $T_0$ and $f$ maps
 every hinge onto a  quadrilateral in $\R^2$.
\end{proposition}

\begin{proof}
Let  $\hat x$ be a point in $\widehat{\Sigma}$. This point belongs
to a simplex $\widehat{\sigma} = (\sigma , [c])$ in
$\widehat{\mathcal{T}}$. Choose an admissible arc $c$ connecting
the base point $x_0 \in T_0$ to $\sigma$. Because $c$ is
admissible, it crosses only finitely many edges $e_1,e_2,...,e_m$
in that order (repetitions may occur). We associate to the path $c$
a sequence of triangles $T_1,T_2,...,T_m \in \mathcal{T}$ by
requiring that $(T_0,T_1 )$ be the hinge of $e_1$, $(T_1,T_2 )$
be the hinge of $e_2$ and so on. We then define $f_j : T_j \to
\R^2$ to be the continuation of $f_{j-1}$  across the edge $e_j$
(for $1 \leq j \leq m$) and we finally set $f(\hat x) =
f_m(P(\hat x)$. The point $f(\hat x)\in \mathbb{R}^2$ only 
depends on the homotopy class $[c]$ and not on the
representative path $c$.

It is clear from the construction that  $f : \widehat{\Sigma} \to
\R^2$ maps every hinge onto a  quadrilateral in $\R^2$. Since $f$
extends $f_0$, the proof is complete.

\end{proof}

\begin{definition}
 The map $f : \widehat{\Sigma} \to \R^2$ is the \emph{development map}\index{development map}
 of the  piecewise flat surface.
\end{definition}

If $f': \widehat{\Sigma} \to \R^2$ is another development, then,
clearly, $f' = g\circ f$ where $g : \R^2 \to \R^2$ is the unique
isometry of the plane such that $g(f(T_0)) = f'(T_0)$.

\medskip

When $\Sigma$ is the boundary surface of a convex polyhedron in
$\R^3$, the development is a very concrete operation. It is
obtained by first placing the initial face (home) somewhere on the
plane and then rolling without slipping the polyhedron, face
after face, following an admissible path. Observe in particular
that we can move our polyhedron toward any point in the plane.
This is a general fact:

\begin{proposition}\label{dev.surjective}
Let $\Sigma$ be a compact piecewise flat surface without boundary. Then any development
$f : \widehat{\Sigma} \to \R^2$ is surjective.
\end{proposition}

\begin{proof} Observe first that $f$ is a closed map (because it is an isometry on each triangle).
Suppose that $\R^2 \setminus f(\widehat{\Sigma}) \neq \emptyset$,
then this set is open and we can find a point $y\in \R^2$ which
lies on the boundary of $f(\widehat{\Sigma})$. Because $f$ is
closed, we can find $\hat x \in \widehat{\Sigma}$ with $f(\hat x)
= y$. Let $x = P(\hat x) \in \Sigma$. This point cannot be in the
interior of any triangle of the triangulation, thus $x$ lies on
an edge $e$. Moving slightly the point $y$ if necessary, we can
assume that $x$ lies in the interior  $e$ (i.e. that $x$ is not a
vertex).

Since  $\Sigma$ has no boundary, $e$ is an interior edge; the
developing map $f$  sends the hinge of $e$ onto a quadrilateral
$Q$ in $\R^2$. The interior of $e$ is sent in the interior of $Q
\subset f(\widehat{\Sigma})$. This contradicts the point $y$
lying on the boundary of  $f(\widehat{\Sigma})$.
\end{proof}

\subsection{The holonomy of a  piecewise flat surface}
\label{sec.holonomy}

The set of all admissible homotopy classes in the piecewise flat
surface $(\Sigma, \mathcal{T})$ which start and end at the base
point $x_0$ form a group $\pi$ with respect to the concatenation. 
This group coincides with the fundamental group $\pi_1(\Sigma',x_0)$.

If $[a]\in \pi$ and $(\sigma, [b]) \in \widetilde{\mathcal{T}}$,
then $(\sigma, [ba])$ is well defined, and this gives a
simplicial action of $\pi$ on $\widetilde{\mathcal{T}}$.

Corresponding to this simplicial action, there is an action of
$\pi$ on $\widehat{\Sigma}$ by isometries; the orbits space of this
action coincides with the surface itself. i.e.  $\Sigma =
\widehat{\Sigma}/\pi$.

\medskip

If  $f : \widehat{\Sigma} \to \R^2$ is a development map  of
$(\Sigma, \mathcal{T})$ and $\gamma = [c]\in \pi$, then there is
a unique isometry $g : \R^2 \to \R^2$ such that the  $g(f(T_0)) =
f(T_0, [c])$. We denote this isometry by $g = \varphi (\gamma)$.

\begin{proposition}
The map $\varphi : \pi \to E(2)$ (the group of isometries of the euclidean plane) is a
group homomorphism.
\end{proposition}

\begin{proof}
 This easily  follows from the construction of the development map.

\end{proof}

\begin{definition}
The homomorphism $\varphi : \pi \to E(2)$ is called the
\emph{holonomy}\index{holonomy} associated to the development $f$.
\end{definition}

\begin{proposition}\label{noboundedorbit}
If $\Sigma$ is compact without boundary, then the group \\ $H =
\varphi (\pi) \subset E(2)$ has no bounded orbit (in particular
it has no fixed point).
\end{proposition}

An obvious but important consequence is the fact that $H$ is not
conjugate to a subgroup of $O(2)$.

\begin{proof}
Suppose that there is a point $y\in \R^2$ such that $H \cdot y$
is bounded. Since the development map is surjective, there exists
a point $\hat x \in \widehat{\Sigma}$ such that $f(\hat x) = y$. 
Observe that $H \cdot y = f(\pi \cdot \hat x)$. Any point in the
surface $\Sigma$ can be connected to $x = P(\hat x)$ by a path of
length at most $D =  \mathrm{diam}(\Sigma)$, hence any point in
$\widehat{\Sigma}$ can be connected to a point in the orbit  $\pi
\cdot \hat x$ by a path of length at most $D$.

Since $f$ preserves the length of all paths, it follows  that any
point in the image $f(\widehat{\Sigma})$ can be connected to a
point in the orbit  $H \cdot y $ by a path of length at most $D$.
The last assertion contradicts the surjectivity of $f$.

\end{proof}

\medskip

Recall that the development of a piecewise flat surface is not
unique, it depends on the choice of an isometry of the home
triangle into $\R^2$. However the holonomy is well defined up to
conjugacy:

\begin{proposition}\label{prop.conjhol}
Let   $f, f': \widehat{\Sigma} \to \R^2$ be two  development
maps  of the piecewise flat surface $(\Sigma, \mathcal{T})$, and
let $\varphi,\varphi' : \pi \to E(2)$ be the corresponding
holonomies. Then $\varphi' (\gamma) = g \varphi(\gamma) g^{-1}$
where $g\in \mathrm{E}(2)$ is the unique isometry such that $f' =
g\circ f$.
\end{proposition}

\begin{proof}
 The holonomy is defined by the condition
 $f (T_0, \gamma) = \varphi(\gamma) (f(T_0))$,
 hence
 \begin{align*}
 \varphi'(\gamma) (f(T_0)) = f' (T_0, \gamma) =
 g\circ f (T_0, \gamma)  &= g\circ \varphi(\gamma) (f(T_0))
 \\ & =  g\circ \varphi(\gamma)\circ g^{-1}(f'(T_0)).
\end{align*}
\end{proof}

\subsection{The development near a singularity}

The previous notions can be clearly visualized if one restricts
one's attention to a simply connected region $\Omega \subset
\Sigma$ which is a union of triangles and which contains exactly
one singular vertex $p$ of order $\b = \b(p) \neq 0$.

Suppose that
the base point $x_0$ sits in $\Omega$ and choose a loop
$c$ in $\Omega' = \Omega \setminus \{p\}$, based at $x_0$ and
surrounding the point $p$ once (so that $[c]$ is a generator of
$\pi_1(\Omega',x_0) \cong \mathbb{Z}$).

Choose a connected component $\widehat{\Omega}$ of the inverse
image $P^{-1}(\Omega) \subset \widehat{\Sigma}$ and still denote
by $P$ the (restriction of the) projection $P : \widehat{\Omega}
\to \Omega$.

We want to describe the geometry of $\widehat{\Omega}$, of the
map $P$ as well as the development and holonomy restricted to
$\widehat{\Omega}$.

It is enough to consider the case where $\Omega$ is the ``star'' of
the vertex $p$, i.e. the union of all triangles incident with $p$
(if $\Omega$ is a larger region, the other triangles will simply
appear as an appendix glued to the star of $p$).

\medskip

The space $\widehat{\Omega}$ is the geometric realization of a
simplicial complex whose simplices are  simplices in $\Omega$
together with an admissible homotopy class of curve connecting
the base point to the given simplex.

\medskip

Let us denote by $T_1,T_2,...,T_k$ the list of all triangles
(i.e. 2-simplices) incident with $p$ and assume that $x_0\in T_1$.
Assume also that $T_i$ has a common edge with $T_{i+1}$ and $T_k$
has a common edge with $T_{1}$. Then a triangle in
$\widehat{\Omega}$ is given by a pair $(T_i,[a])$ where $[a]$ is
the homotopy class of a curve $a$ in $\Omega'$ from $x_0$ to
$T_i$. This homotopy class is parametrized by a single integer
$d\in \mathbb{Z}$ (the degree of $a$) which counts the number of
times $a$ turns around the point $p$. In other words,
$\widehat{\Omega}$ is an infinite strip made out of  countably
many copies of each triangle  $T_1,T_2,...,T_k$ each indexed by
the degree $d\in \mathbb{Z}$
$$
 \widehat{\Omega} = \bigcup_{d\in \mathbb{Z}}
\left( T_{1,d} \cup T_{2,d} \cup ... \cup T_{k,d}\right).
$$
To develop $ \widehat{\Omega}$, start with an isometry $f_1$ from
$T_1$ to a triangle in the euclidean plane and continue this
isometry by unfolding each hinge in $ \widehat{\Omega}$. The
development $f$ then clearly satisfies
$$
 f(T_{i,d}) = R^d f(T_{i,0})
$$
where $R$ is a rotation of angle $\theta$ ($=$ the sum of the angles at $p$ of the triangles  $T_1,T_2,...,T_k$) around the point $q =
f_1(p)$. The rotation $R\in E(2)$ is clearly the holonomy of the
generator $[c]$ of $\pi_1(\Omega',x_0)$.

\medskip

We collect in the next proposition, some of the conclusions of the previous discussion:

\begin{proposition}\label{local.dev}
\begin{enumerate}
  \item The inverse  image $P^{-1}(p)$ of $p$ in $\widehat{\Omega}$ contains exactly one point $\hat p$;
  \item the holonomy   $\varphi (c)$ of $[c]$  is a rotation of angle $\theta = 2\pi(\b +1)$;
  \item  if $\b$ is not an integer, then $q=f(\hat p)$ is the unique fixed point of the rotation  $\varphi (c)$.
\end{enumerate}
\end{proposition}

 \qed


\subsection{Geometric equivalence of euclidean triangulations}

Let  $(\Sigma,\Tau)$ be a piecewise flat surface. Choose a triangle $T_{\alpha_0} \in \Tau$ and a point $q$ in the interior of an edge of $T_{\alpha_0}$. If one connects the point $q$ to the opposite vertex in $T_{\alpha_0}$ by a euclidean segment, one obtains
two subtriangles $T'_{\alpha_0}, T''_{\alpha_0}$ whose union
is $T_{\alpha_0}$.

If one replaces  the triangle $T_{\alpha_0}$ with $T'_{\alpha_0}$ and $T''_{\alpha_0}$ in the  triangulation  $\Tau_{q}$, one obtains  a new  triangulation  $\Tau_{q}$.

\begin{definition}
a) The triangulation  $\Tau_{q}$  is said to be obtained from  $\Tau$ by an \emph{elementary subdivision}.

b) The \emph{geometric equivalence} is the equivalence relation on the set
of euclidean triangulations on a surface which is generated by
elementary subdivisions.
\end{definition}

In other words, two  euclidean triangulations $\mathcal{T}_1 ,\mathcal{T}_2$
on the surface $\Sigma$ are geometrically  equivalent if
there is a common subdivison $\mathcal{T}$ which is also
a euclidean triangulation.

\begin{proposition}\label{prop.geom.inv}
The area measure $dA$, the length structure $\ell$,  the singularity order $\beta$, the development and holonomy are invariants of this equivalence relation.
\end{proposition}

\begin{proof}
The statement is obvious for  $dA$,  $\ell$ and $\beta$.
Observe now that if $T$ is a triangle of $\Tau$ and
 $T', T''$ is an elementary subdivision of $T$, then
the pair $(T',T'')$ is the hinge of their common edge $e \subset
T$. Observe also that if   $f : T \to \R^2$ is an isometry, then
$f$ is an unfolding of that hinge.

This argument shows that
the development remains unchanged when subdivising the triangulation.
Since the development is invariant, so is the holonomy.

\end{proof}


\subsection{Flat metrics with conical singularities}

If $(\Sigma, \mathcal{T})$ is a piecewise flat surface, then $\Sigma'$ carries a well defined riemannian metric $m$; this metric is flat (i.e. it has no curvature) and in  the neighbourhood of a conical singularity of total angle $\theta$, we can introduce polar coordinates $(r,\varphi)$, where $r\geq 0$ is the
distance to $p$ and  $\varphi\in \R/(\theta \mathbb{Z})$ is the angular variable (it is defined modulo $\theta$).
In these coordinates, the metric reads
$$
 m = dr^2 + r^2d\varphi^2.
$$
A calculation shows that this metric can be written as
\begin{equation}\label{ds2beta}
  m_{\b} = |z|^{2\beta}|dz|^2,
\end{equation}
where $ z = \frac{1}{\b+1}(r e^{i\varphi})^{\b+1}$ (see \cite{Tr86}).

\begin{definition}
A \emph{flat surface with conical singularities}  $(\Sigma, m)$ is a surface $\Sigma$ together with
a singular Riemannian metric $m$ which is isometric
to the metric $m_{\b}$ in (\ref{ds2beta}) in the
neighbourhood of every point $p\in \Sigma$, where $\b = \b(p) \in (-1,\infty)$.

\medskip

One says that $\b(p)$ is the \emph{singularity order}\index{singularity order}
of $p$ and  $p$ is a conical  singularity if $\b(p) \neq 0$.
The singular points form a discrete set and the formal sum (with discrete support)  $\sum\b(p)\; p$ is called the \emph{divisor} of the singular
metric $m$. One sometimes also says that $m$ \emph{represents}  this divisor.
\end{definition}

\begin{proposition}\label{triangulable}
Any compact flat surface with conical singularities $(\Sigma, m)$ can be geodesically triangulated. The resulting triangulation is
a euclidean triangulation on $\Sigma$ and the associated length structure coincides with the length in the metric $m$.
\end{proposition}

A proof can be found in \cite{Tr86} and in \cite{thurston98}. See also \cite{ILCT}
and \cite{rivin} for further discussions on triangulations of piecewise flat surfaces.

\qed

\begin{proposition}\label{metrisable}
Two euclidean triangulations $\Tau,\Tau'$ on a compact surface $\Sigma$ are geometrically equivalent if and only if
they give rise to the same flat surface with conical singularities $m$ on $\Sigma$.
\end{proposition}

\emph{Proof.}  It is clear from  Proposition \ref{prop.geom.inv}, that two triangulations which are geometrically equivalent give rise to the same singular  flat metric. Conversely, suppose that the triangulations $\Tau$ and $\Tau'$ define the same metric, then each triangle of $\Tau$ is decomposed by $\Tau'$
in a finite number of polygonal regions. We may then further decompose these polygons in euclidean triangles, and we thus obtain a new euclidean triangulation of $\Sigma$ which is a subdivision of both $\Tau$ and $\Tau'$.

\qed

\subsection{Relation with Riemann surfaces}

If  $(\Sigma, m)$ is an oriented  flat surface with conical singularities, then it is covered by charts $\{(U_j,z_j)\}$ such
that the metric $m$ takes the form (\ref{ds2beta}) in each $U_j$.
The transition from one such  coordinate $z_j$ to another one is given by a conformal transformation. Thus $\Sigma$ is  a Riemann surface with a holomorphic atlas given by $\{(U_j,z_j)\}$.

\begin{remark}\label{sing_supprimable}
The reader should observe here that the conical singularities are invisible from the conformal viewpoint. This is 
a consequence of the formula  (\ref{ds2beta}) which shows that the singular metric is conformal to a smooth metric. 
It can also be seen as a consequence of the theorem of
removability of singularities of locally bounded meromorphic functions. 
\end{remark}

In the converse direction, we can start with a closed Riemann surface with a divisor and ask whether there is a conformal
flat metric representing this divisor. The answer is positive and
the following theorem  classifies all compact euclidean
surfaces with conical singularities.

\begin{theorem}\label{th.srsc}
Let $\Sigma$ be a compact connected Riemann
surface without boundary. Fix $n$ distinct points $p_1,p_2, \dots , p_n \in \Sigma$ and
$n$ real numbers $\b_1,\b_2,\dots ,\b_n \in
(-1,\infty )$.

There exists a conformal flat metric $m$ on $\Sigma$ having a conical
singularity of order $\b_j$ at $p_j$ ($j = 1,\dots , n$) if and only if
the Gauss-Bonnet condition $\chi(S) + \sum_{j=1}^n\b_j = 0$  holds.
This metric is unique up to homothety.
\end{theorem}

See  \cite{Tr86}, a shorter proof can be found in  \cite[\S IV]{Tr93}.

\qed

\begin{remark}
A careful examination of the proof shows that the metric $m$
depends continuously on all the parameters: The conformal structure, the points $p_j$ and the orders $\b_j$.
\end{remark}

{\footnotesize
There is a similar  theorem  for the case of hyperbolic
metrics with conical singularities, see
\cite{Heins,McO,Picard,Tr91}. There are  are also
various other extensions (non constant curvature, non orientable surfaces, non compact surfaces, and surfaces with boundary, see \cite{ht92,Tr91, Tr93}). The case of spherical metric is more delicate, see \cite{eremenko, uy} for a
study of spherical metric with three conical singularities on
the 2-sphere.
}

\medskip

\begin{theorem}
Given a compact oriented surface $\Sigma$, there are natural  bijections
between the following three sets:
\begin{enumerate}[1)]
  \item The set of geometric equivalence classes of euclidean triangulations   on $\Sigma$ up to homothety;
  \item the set of flat metrics $m$ on $\Sigma$ with conical
singularities up to homothety;
  \item the set of conformal structures on $\Sigma$ together
  with a finite real  divisor $\sum_i \b_ip_i$ such that $\b_i > -1$ and
  the Gauss-Bonnet condition (\ref{cond.GB}) is satisfied.
\end{enumerate}
\end{theorem}

\begin{proof}
Theorem \ref{th.srsc} says precisely that there is a bijection between sets (2) and (3). Proposition \ref{metrisable}
shows that there is a natural injection from (1) to (2),
this injection is surjective by Proposition
 \ref{triangulable}.

\end{proof}

\section{Punctured surfaces}

\subsection{Punctured surfaces and their fundamental groups}

\begin{definition}
We define a \emph{punctured surface}\index{punctured surface}  $\Sigma_{g,n}$ to be an oriented, closed connected surface $\Sigma$  of genus $g$ together
with a distinguished set of $n$  pairwise distinct points $p_1,p_2,...,p_n \in \Sigma_{g,n}$.
\end{definition}
The points $p_1,...,p_n$ are considered to be special places (with some geometric significance) on the surface. We call them the  \emph{punctures} and we  denote by $\Sigma'_{g,n}$ the surface obtained by removing them:
$$
 \Sigma'_{g,n} = \Sigma_{g,n}\setminus \{ p_1,p_2,...,p_n \}.
$$
The connected sum of two punctured surfaces is defined by removing a disk
containing no puncture in each surface and then gluing them along their
boundary. The resulting surface is again a punctured surface. In fact we have
$$
 \Sigma_{g_1,n_1} \, \# \, \Sigma_{g_2,n_2} = \Sigma_{g_1+g_2,n_1+n_2},
$$
where the symbol $\#$  means the connected sum. 
In particular
\begin{equation}\label{dec.Psurface}
\Sigma_{g,n} = \Sigma_{g,0} \, \# \,\Sigma_{0,n}.
\end{equation}

We easily deduce from this observation  that the Euler characteristic of $\Sigma'_{g,n}$ is given by
\begin{equation}\label{ }
 \chi (\Sigma'_{g,n}) = 2-2g-n.
\end{equation}
If $n>0$, then $\Sigma'_{g,n}$ can be homotopically retracted onto a bouquet of  $2g+n-1$ circles and
 the  fundamental group $\pi_{g,n}$  of $\Sigma'_{g,n}$
is thus a free group on $2g+n-1$ generators.

Note  that $\pi_{g,n}$  also admits the following presentation with
$2g+n$ generators and one relation:
\begin{equation}\label{can.presentation}
\pi_{g,n} = \langle
 a_1,...,a_g, b_1,...,b_g, c_1,...,c_n   \tq
\Pi \left[a_i,b_i\right] = \Pi c_j
\rangle,
\end{equation}
this presentation is a consequence of the identity (\ref{dec.Psurface})  and Van Kampen's Theorem.

\subsection{Uniformization of a punctured Riemann surface}
\label{subs.uniformization}

Let us fix a conformal structure
$[m]$ on $\Sigma_{g,n}$. Assuming that $2-2g-n<0$,  the  uniformization Theorem states that $(\Sigma', m)$ is  conformally equivalent to $\mathbb{U}/\Gamma$ where
$\mathbb{U} = \{ z \in \mathbb{C} \tq \mathrm{Re}(z)>0\}$ is the upper-half plane, and $\Gamma\subset PSL_2(\R)$ is a Fuchsian group of the first kind\footnote{Recall that a \emph{Fuchsian group} is a discrete subgroup of $PSL_2(\R)$, it is of the fi\emph{rst kind} if there is a fundamental domain $D\subset \mathbb{U}$ of finite hyperbolic area.} isomorphic to $\pi_{g,n}$.

\medskip

The isomorphism $\pi_{g,n} \to \Gamma$ is compatible with the punctures  in the sense that  the generator $c_i$ is sent to a parabolic element of $\Gamma$ and the generators $a_i,b_i$ are sent to  hyperbolic elements  (here, the letters $a_i,b_i,c_j$ refer to the presentation (\ref{can.presentation})).

\medskip 

Let us denote by  $\Upsilon \subset \partial \mathbb{U} = \R \cup \{ \infty\}$ the
set \emph{cusp points}\index{cusp points}  of $\Gamma$, i.e. the set of  fixed points of all parabolic elements in $\Gamma$. Following  \cite[page 10]{shimura}, we define a topology on $\widehat{\mathbb{U}} = \mathbb{U}\cup \Upsilon$
as follows: for a point $z\in \mathbb{U}$, the family of hyperbolic disks $D(z,\rho)$ is a fundamental system of neighborhoods of $z$. For  a point $y\in \Upsilon$ the family of horodisks centered at $y$ is a fundamental system of neighborhoods of $y$. With this topology, $\widehat{\mathbb{U}}$ is a Hausdorff space and $\Gamma$ acts by homeomorphisms. The space is not locally compact and $\Upsilon$ is
topologically a discrete space.  Standard arguments from hyperbolic geometry (see e.g. \cite{shimura})  show that the projection map $P : \mathbb{U} \to  \mathbb{U} /\Gamma = \Sigma'$ extends to a surjective continuous map
\begin{equation}\label{Udevmap}
 \widehat{P} :  \widehat{\mathbb{U}} \to \Sigma
\end{equation}
where $\widehat{\mathbb{U}} = \mathbb{U}\cup \Upsilon$. This extension maps $\Upsilon$ to the punctures
$\{p_i\} \subset \Sigma$.

\bigskip

\textbf{Remarks} \textbf{1)} The previous considerations show that
there exists a unique  metric $m_{-1} $ on $\Sigma'$
of constant curvature $-1$ which is complete, has finite volume and belongs to the conformal structure $[m]$. 
This metric has a cusp at each puncture $p_i$, it is the unique metric such that
$P^*m_{-1}$ is the Poincar\'e metric on $\mathbb{U}$; its existence can also be proved by  directly solving the prescribed curvature equation (see \cite{ht90,ht92}).

\medskip

\textbf{2)} We know from Theorem \ref{th.srsc} that
the conformal class $[m]$ also contains
a metric  $m_0$ on $\Sigma'$, unique  up to homothety, which is flat and has a conical singularity of order $\b_j$ at $p_j$ ($j = 1,\dots , n$) provided (\ref{cond.GB}) holds.

This flat metric lifts as a flat conformal metric $\widetilde{m}_0
= P^*(m_0)$ on $\mathbb{U}$. For this metric, $\mathbb{U}$ is not complete and its completion is given by $\widehat{\mathbb{U}}=\mathbb{U}$. The map
$\widehat{P} :  \widehat{\mathbb{U}} \to \Sigma$ is thus a
concrete model of the universal branched covering introduced
earlier.

\medskip

We identify the set $\Upsilon$ as a subset of $\Gamma$ as follows: we first
fix a base point $\tilde{z}_0 \in \mathbb{U}$ and let $z_0 = P(\tilde{z}_0) \in \Sigma$. For $y \in \Upsilon$, let us denote by $\tilde{\gamma}_y$ the hyperbolic ray in $\mathbb{U}$ starting at  $\tilde{z}_0$ and asymptotic to the point $y$, and let ${\gamma}_y = P(\tilde{\gamma}_y)$, this is a path joining
$z_0$ to a puncture  $p_i = P(y)$. Now let $D_i \subset \Sigma$ be a small disk around $p_i$ containing no other puncture, and let $\gamma'_y = \gamma_y \setminus D_i$.

We now define $c_y\in \pi_1(\Sigma', z_0)$ to be the homotopy class of the path obtained by following $\gamma'_y$, then $\partial D_i$ (in the positive direction) and then
$(\gamma'_y)^{-1}$.

Recall that we have a canonical isomorphism,  $\Gamma \cong \pi_1(\Sigma', z_0) = \pi_{g,n}$, we have thus constructed a map
\begin{equation}\label{Ydevmap}
\left.\begin{array}{ccc}\Upsilon & \to & \Gamma \\
y & \mapsto & c_y\end{array}\right.
\end{equation}
It is clear that $c_y \in \Gamma$ is a parabolic element fixing $y$, in particular, the map $\Upsilon \to \Gamma$ is injective.

\subsection{Some  Groups of Diffeomorphisms of a Punctured surface}

Given a punctured surface $\Sigma = \Sigma_{g,n}$, we define
$\Diff_{g,n}$ to be the group of diffeomomorphisms $h : \Sigma \to \Sigma$ which leaves the set $\{p_1,...,p_n\}$ of punctures  invariant. We also introduce the following subgroups: $\Diff^+_{g,n}\subset \Diff_{g,n}$ is the subgroup of orientation preserving diffeomomorphisms, $\PDiff_{g,n}$ is the subgroup of pure diffeomomorphisms, i.e. diffeomomorphisms fixing each puncture $p_i$ individually and
$\PDiff^+_{g,n}= \PDiff_{g,n}\cap \Diff^+_{g,n}$.

\medskip

Every element $h\in \Diff^+_{g,n}$ permutes the punctures
and we have an exact sequence
$$1 \to \PDiff^+_{g,n} \to \Diff^+_{g,n} \to \Sym(n) \to 1.$$
where $\Sym(n)$ is the  permutation group of
$\{p_1,...,p_n\}$.

We also define $\Diff^0_{g,n}\subset \PDiff^+_{g,n} $ to be the group of diffeomorphisms which
are isotopic to the identity through an
isotopy fixing the punctures. The quotient
$$
 \Mod_{g,n} = \pi_0(\Diff^+_{g,n}) = \Diff^+_{g,n}/\Diff_{g,n}^0,
$$
is called the \emph{mapping class group}\index{mapping class group} or the 
\emph{modular group}\index{modular group}
of the punctured surface $\Sigma_{g,n}$, and
$$
 \PMod_{g,n} = \pi_0(\PDiff^+_{g,n}) = \PDiff^+_{g,n}/\Diff_{g,n}^0,
$$
is the \emph{pure mapping class group}\index{pure mapping class group}.

\medskip

These groups have been intensely studied since the pioneer work of Dehn and Nielsen.
We refer to \cite{birman74,ivanov02,morita,zvc} among many other papers,  for more information.

\subsection{Outer automorphisms}

The mapping class group is related to the group of outer automorphisms of
the fundamental group of $\Sigma'$. Let us recall this purely algebraic notion: If  $\pi$ is an arbitrary group, we denote by $\Aut (\pi)$  the group of all its automorphisms and
by $\Inn (\pi)\subset\Aut(\pi)$ the subgroup of inner automorphisms (i.e. conjugations $\gamma \to \alpha \gamma \alpha^{-1}$).  This is a normal subgroup.

\begin{definition}
The group of   \emph{outer automorphisms} \index{outer automorphisms}   of $\pi$  is the quotient
$$\Out(\pi)  =  \Aut(\pi)/\Inn(\pi).$$
An outer automorphism is thus an automorphism of $\pi$  defined up to conjugacy.
\end{definition}

\begin{lemma}\label{lem.DNB}
There is a naturally defined group homomorphism
 $$
 \Mod_{g,n} \to \Out(\pi_{g,n}).
 $$
\end{lemma}

\begin{proof}
This homomorphism is defined as follows. Let $h \in \Diff(\Sigma')$
be an arbitrary diffeomorphism and fix a  base point $*$
and a path $\delta$  in  $\Sigma'$ connecting $*$ to $h(*)$.
If $\gamma$ is a loop in $\Sigma'$ based at $*$, then we set
$$h_{\delta} (\gamma) = \delta^{-1}  (h \circ \gamma) \delta.$$
This defines an automorphism $h_{\delta, \#} \in \Aut(\pi_{g,n})$.

If  $\delta'$  is another  path connecting $*$ to $h(*)$, then
$h_{\delta, \#}$ and $h_{\delta', \#}$ are conjugate by $\delta^{-1}\delta'$.
The outer automorphism $h_{\#} \in \Out(\pi_{g,n})$
is thus well defined independently of the choice of the path $\delta$ and it  
is clear that if  $h$ is homotopic to the identity, then it  acts trivially on $\pi_{g,n}$, i.e. we have defined a map
$\Mod_{g,n} \to \Out(\pi_{g,n})$. It is routine to check  that
it is a group homomorphism.
\end{proof}

\medskip

Introducing the group  $\POut(\pi_{g,n}) \subset \Out(\pi_{g,n})$ of all
outer automorphisms preserving the conjugacy class of each $c_i$ 
($i=1,...,n$) in the presentation (\ref{can.presentation}), we have the following deep result:

\begin{theorem}\label{th.DNB}
If $g>0$ and $n>0$, then the  homomorphism defined in the previous lemma
 induces an isomorphism
\begin{equation}\label{eq.DNB}
  \Phi : \PMod_{g,n} \stackrel{\sim}{\longrightarrow} \POut(\pi_{g,n}).
\end{equation}
\end{theorem}

This is the so called \emph{Dehn-Nielsen-Baer Theorem}, see
\cite{ivanov02,zvc} for a proof.

\qed

\subsection{Lifting the group $\Diff^0(\Sigma')$  on $\mathbb{U}$ }

Using the notations of section \ref{subs.uniformization},
one writes the universal branched covering of  
$\Sigma_{g,n}$ as $\widehat{P} :  \widehat{\mathbb{U}} \to \Sigma$, where $\widehat{\mathbb{U}}=\mathbb{U}\cup \Upsilon$
(we still assume  $2-2g-n<0$).

We denote by $\Diff^+(\mathbb{U})$ the group of  orientation preserving diffeomorphisms
of $\mathbb{U}=\widetilde{\Sigma'}$ and we define the \emph{normalize}r $N(\Gamma)$
and the \emph{centralizer} $C(\Gamma)$ of $\Gamma$  in $\Diff^+(\mathbb{U})$ by
$$
 N(\Gamma) = \{ h \in \Diff^+(\mathbb{U}) \tq h\Gamma = \Gamma h \}.
$$
and
$$
 C(\Gamma) = \{ h \in \Diff^+(\mathbb{U}) \tq h\circ \gamma = \gamma \circ h \mbox{ for all } \gamma \in \Gamma\}.
$$
Observe that $C(\Gamma) = \ker (\psi)$, where $\psi : N(\Gamma) \to \Aut(\Gamma)$ is defined by $\psi(h) : \gamma \to h \gamma h^{-1}$.

The \emph{center} of $\Gamma$ is the intersection $Z(\Gamma) = \Gamma\cap C(\Gamma)$; it is the largest abelian subgroup of $\Gamma$.

\begin{lemma}
 Let $\Gamma$ be an arbitrary Fuchsian group, then $Z(\Gamma)$ is trivial unless $\Gamma$ is cyclic.
\end{lemma}

\begin{proof} This follows from classical  Fuchsian group theory. Indeed, it is well known that if  $\gamma_1,\gamma_2 $ are non-trivial elements in  $PSL_2(\R)$, then they commute  if and only if they have the same fixed points (see e.g.  \cite[theorem 2.3.2]{katok}). So if $Z(\Gamma)$ contains a non-trivial element $\gamma_0$, then  any $\gamma \in \Gamma\setminus \{ id\}$ must have the same
fixed points as $\gamma_0$ and it follows from \cite[theorem 2.3.5]{katok}) that $\Gamma$ is cyclic.

\end{proof}

\medskip

Recall the projection  $P :\mathbb{U} \to \Sigma' = \mathbb{U}/\Gamma$.
For any element $h\in N(\Gamma)$,  we define $P_*h : \Sigma' \to \Sigma' $
by $P_*h(x) = P(h(\tilde{x}))$ where $\tilde x \in \mathbb{U}$ is an arbitrary point in   $P^{-1}(x)$. This map is well-defined, because the condition  $h\Gamma = \Gamma h$ means precisely that $h$ maps $\Gamma$-orbits in $U$ to $\Gamma$-orbits, and it is clearly a diffeomorphism.
We thus have defined a map
$$P_* : N(\Gamma) \to \Diff^+(\Sigma'),$$
and it is obviously a group homomorphism.

\begin{proposition}
$\Gamma$ is a normal subgroup in $N(\Gamma)$ and $P_*$
defines an isomorphism from
$N(\Gamma)/\Gamma$  to $\Diff^+(\Sigma')$.
\end{proposition}

\begin{proof}
It is obvious that $\Gamma\subset N(\Gamma)$ is normal and that
$P_*(\Gamma) = \{id\}$. In particular $P_*$ factors through a well
defined  homomorphism  $N(\Gamma)/\Gamma\to \Diff^+(\Sigma')$.
This homomorphism is surjective since every diffeomorphism of
$\Sigma'$ lifts to the universal cover $\mathbb{U}$ of $\Sigma'$.

Suppose now  that $P_*h = id$. Then $h(x) \in \Gamma
\cdot x$ for all $x\in U$. This means that there exists a map $U
\to \Gamma$, $x \to \gamma_x$ such that  $h(x) =  \gamma_x\, x$
for all $x\in U$. Since $h$ is continuous, so is this map, but
this implies that $x \mapsto \gamma_x$ is constant because $\Gamma$
is a discrete group. It follows that $h \in  \Gamma$ and we have
shown that $P_*  : N(\Gamma)/\Gamma\to \Diff^+(\Sigma')$ is also
injective.

\end{proof}

\begin{lemma}
 $P_*$ maps $C(\Gamma)$   isomorphically onto $\Diff^0 
 (\Sigma')$.
\end{lemma}

\begin{proof}
Suppose that  $P_*h \in \Diff^0 (\Sigma')$. Then there exists an isotopy
$h_t \in N(\Gamma)$ such that $h_0 = id$ and $h_1=h$. Hence $\psi(h_t) \in \Aut(\Gamma)$ is constant by continuity. Because $\psi(h_0) = \psi(id)\in \Aut(\Gamma)$ is the trivial element, we have $h \in \ker \psi = C(\Gamma)$.

 \medskip

In the reverse direction, we use an argument going back to
Nielsen: Suppose that  $h \in \ker \psi = C(\Gamma)$ and define
$h_t(x)\in U$ to be the point on the hyperbolic segment
$[x,h(x)]$  such that $d(x,h_t(x)) = t d(x,h(x))$ (where $d$ is
the hyperbolic distance in $\mathbb{U}$). Since $h \in \ker \psi$
and $\Gamma$  preserves the hyperbolic distance in $\mathbb{U}$,
the segment $[\gamma x,h(\gamma x)]$ coincides with $[\gamma x,\gamma
h(x)]$ for any $x\in \mathbb{U}$   and any $\gamma \in \Gamma$.
Therefore we have  $h_t (\gamma x) = \gamma h_t (x)$, i.e. $h_t
\in C(\Gamma)\subset N(\Gamma)$. The path $P_*h_t \in
\Diff(\Sigma')$ is an isotopy from $P_*h$ to the identity and we
conclude that $P_*h \in  \Diff^0(\Sigma')$.

\medskip

We have proved that $P_*^{-1}(\Diff^0 (\Sigma')) = C(\Gamma)$.
It is now clear that $P_*: C(\Gamma) \to \Diff^0 (\Sigma')$ is an isomorphism since its kernel is   $C(\Gamma)\cap \Gamma = Z(\Gamma) =\{id\}$.

\end{proof}

\begin{corollary}
$P_*$  induces an isomorphism from  $N(\Gamma)/(\Gamma \times C(\Gamma))$ to the modular group $\Mod_{g,n}$.
\end{corollary}


\section{The representation variety of a finitely generated group in $SE(2)$} \label{sec.rep}

Given a finitely generated group $\pi$ and an algebraic Lie group
$G$, it is easy to see that the set $\Hom(\pi,G)$ is an algebraic set. The group $G$ itself acts on $\Hom(\pi,G)$ by conjugation:
$g\cdot \varphi (\gamma) = g^{-1} \varphi (\gamma) g$. The quotient  space is called the \index{representation variety} \emph{representation variety} of $\pi$ in
$G$ and denoted by 
$$\mathcal{R}(\pi  , G)=\Hom(\pi,G)/G.$$
This variety plays an important role in the study of geometric structures on manifolds, see e.g. \cite{goldman88}.

\medskip

The discussion in section \ref{sec.holonomy} shows that an
element of the representation variety 
$\mathcal{R}(\pi  , \mathrm{E}(2))$ is associated to any
 piecewise flat surface $(\Sigma, \mathcal{T})$
 (where $\pi = \pi_1(\Sigma',x_0)$). In the present section, we
 investigate the structure of $\mathcal{R}(\pi  , \mathrm{E}(2))$, 
 in fact, for  convenience, we shall restrict ourself to the subgroup
$\SE(2) \subset \mathrm{E}(2)$ of orientation preserving isometries of the euclidean plane (this is a subgroup of index 2).

\subsection{On the cohomology of groups}

We will need some elementary results from group cohomology;
here we recall a few basic definitions and facts.

\medskip

Let $\pi$ be an arbitrary group and $A$ be a $\pi$-module, i.e. an abelian group with a representation $\rho : \pi \to \Aut (A)$.

\begin{definition}
 \begin{enumerate}
  \item A \emph{$1$-cocycle}\index{cocycle}   in $A$ is a map
  $\sigma : \pi \to A$ such that
  $$\sigma (\gamma_1\gamma_2) = \sigma  (\gamma_1) +
   \rho(\gamma_1)\cdot\sigma (\gamma_2)
  $$
  for any $\gamma_1,\gamma_2 \in \pi$. The set of $1$-cocycles is an abelian
  group denoted by $Z^1(\pi, A)$.
  \item The $1$-cocycle $\sigma\in Z^1(\pi, A)$ is a \emph{$1$-  coboundary}\index{coboundary}  if it can be written as
  $$\sigma = \delta_a (\gamma) = \rho(\gamma) \cdot a - a.$$
for some  element $a\in A$.   The set of $1$-coboundaries is a subgroup of $Z^1(\pi, A)$ denoted by   $B^1(\pi, A)$.

\item The quotient
  $$H^1(\pi, A) = Z^1(\pi, A)/B^1(\pi, A) $$
is the \emph{first cohomology group\index{group cohomology} of $\pi$ with values in $A$.}
\end{enumerate}
\end{definition}

\medskip

\textbf{Example.}  Let us compute the first cohomology group when
$A=k$ is a field and $\pi$ is a finitely generated group.
We denote by $k_{\rho}$  the $\pi$-module $k$  with the representation $\rho : \pi \to \Aut (k)$.

\medskip

Assume first that the representation $\rho : \pi \to \Aut (k)$ is a scalar representation, i.e. $\rho : \pi \to  k^* \subset \Aut (k)$ and that
$\pi = F_{s} = \langle  a_1, a_2, ..., a_{s}  \rangle$ is a free group on $s$ generators.

\medskip

Since $\pi$ is free, the  homomorphism $\rho : \pi \to k^*$ is completely determined by the vector $r = ( \rho(a_1),\rho (a_2), ..., \rho(a_{s})) \in (k^*)^{s}$. Likewise, a cocycle is given by the vector
$$t = ( \tau(a_1),\tau (a_2), ..., \tau(a_{s})) \in k^{s}.$$
There is no restriction on the vector $t\in k^{s}$ (again because $\pi$ is free) and thus
\begin{equation}\label{eq.Z1}
 Z^1(\pi, k_{\rho}) \cong   k^{s}.
\end{equation}

An element $\sigma \in Z^1(\pi, k_{\rho})$ is a coboundary if $\sigma = u(id-\rho)$ for some $u\in k$, thus
$$
 B^1(\pi, k_{\rho}) \cong  k \cdot  (1-\rho(a_1), 1-\rho (a_2), ...,1- \rho(a_{s}))   \in k^{s}.
$$

Let us choose a  linear form $\mu : k^{s} \to k$ such that $\mu \equiv 0$
if $\rho$ is trivial and
$$\mu (1-\rho(a_1), 1-\rho (a_2), ...,1- \rho(a_{s}))  \neq 0$$
else.
It is easy to check that
$$
B^1(\pi, k_{\rho}) \oplus
\ker \mu = k^{s} = Z^1(\pi, k_{\rho})
$$
in $k^{s}$  and we thus obtain the following
\begin{proposition}\label{compute.H1}
For any free group on $s$ generators, we have
$$
H^1(\pi, k_{\rho}) = Z^1(\pi, k_{\rho})/B^1(\pi, k_{\rho})  = \ker \mu
\cong\begin{cases}
k^{s} &  \quad \text{if $\rho$ is trivial,} \\
k^{s-1} & \quad \text{otherwise. }
\end{cases}
$$
\end{proposition}
\qed

\medskip

\textbf{General case.}  Let us compute the first cohomology group when
$A=k$ is a field and $\pi$ is a finitely presented group with
presentation
$$
\pi = \langle  S  \tq R  \rangle.
$$
Here $S = \left\{a_1, a_2, ..., a_{s} \right\}  \subset \pi$ is a finite set generating the group and  $R = \left\{r_1, r_2, ..., r_m  \right\}  \subset  F(S)$ (= the free group on $S$) is a finite set of words in $S$ defining all the relations among the elements of $S$.
We denote by $k_{\rho}$  the $\pi$-module $k$  with the representation $\rho : \pi \to \Aut (k)$.

\medskip

For any relation $r=a_{i_1}a_{i_2}\cdots a_{i_p} \in R$, we introduce the linear
form $\lambda_r : k^s \to k$ defined by
\begin{equation}\label{ }
 \lambda_r (\sigma)= \sum_{\mu=1}^p \left(\prod_{\nu < \mu} \rho(a_{i_{\nu}}) \right)\sigma(a_{i_{\mu}})
\end{equation}
and we define $\Lambda : k^s \to k^m$, by
\begin{equation}\label{def.Lambda}
  \Lambda (\sigma) = ( \lambda_{r_1} (\sigma), ...,  \lambda_{r_m} (\sigma)).
\end{equation}

\begin{lemma} Th space of $1$-cocycles in $k_{\rho}$ is given by
 $$
  Z^1(\pi,k_{\rho}) = \ker \Lambda =  \bigcap_{r\in R} \ker \lambda_r \subset k^s.
 $$
\end{lemma}

\begin{proof}
If $\sigma \in Z^1(\pi,k_{\rho})$ and $r=a_{i_1}a_{i_2}\cdots a_{i_p} \in R$, then we deduce from the cocycle relation that
\begin{align*}\label{}
   0 & = \sigma (r) =  \sigma(a_{i_1}a_{i_2}\cdots a_{i_p} )  =
   \sigma(a_{i_1}) +  \rho(a_{i_1})\sigma(a_{i_2}\cdots a_{i_p} ) \\
    &  =  \sigma(a_{i_1}) +  \rho(a_{i_1})\sigma(a_{i_2})
    +  \rho(a_{i_{i_1}})\rho(a_{i_2})\sigma(a_{i_3}\cdots a_{i_p} ) \\
    & = \sum_{\mu=1}^p \left(\prod_{\nu < \mu} \rho(a_{i_{\nu}}) \right)\sigma(a_{i_{\mu}}).
\end{align*}
\end{proof}
On the other hand, since any  $1$-coboundary  in $k_\rho$ is a multiple of
$\rho - 1$, we have
$$
   B^1(\pi,k_{\rho}) = k \cdot (\rho - 1) \subset k^s.
$$
We have  proved the following

\begin{proposition}
 The first cohomology group of the finitely presented group
  $\pi =  \langle  S  \tq R  \rangle$ with value in $k_{\rho}$ is given by
  $$
  H^1(\pi,k_{\rho}) =  \ker \Lambda /(k(\rho - 1)).
  $$
\end{proposition}
\qed

In particular, if $\pi$ has exactly one non trivial relation, then
$$
H^1(\pi, k)  \cong\begin{cases}
k^{s-1} &  \quad \text{if $\rho$ is trivial,} \\
k^{s-2} & \quad \text{otherwise. }
\end{cases}
$$
where $s = \mathrm{Card} (S)$ is the number of generators.

\subsection{Abelian Representations}

Representations of a finitely presented group $\pi$ in an abelian Lie group are easy to describe:

\begin{lemma}
If $G$ is an abelian group, then $\mathcal{R}(\pi, G) =
\Hom(\pi, G)$. This set is itself an abelian topological group.
\end{lemma}

\begin{proof}
Since there are no non trivial inner automorphisms in an abelian group,
it is clear that $\mathcal{R}(\pi, G) = \Hom(\pi, G)$.

We endow $\Hom(\pi, G)$ with the compact open topology and we define a product on this space by
$$
  (\varphi_1\varphi_2) (\gamma) =  \varphi_1 (\gamma) \varphi_2   (\gamma)
$$
for $\varphi_1,\varphi_2 \in  \Hom(\pi, G)$ and $\gamma \in \pi$.
The following calculation shows that $\Hom(\pi, G)$ is a group for
this multiplication:
\begin{eqnarray*}
  (\varphi_1\varphi_2) (\gamma_1\gamma_2 ) & = &
\varphi_1(\gamma_1\gamma_2 ) \varphi_2 (\gamma_1\gamma_2 )
\\  & = &
\varphi_1(\gamma_1) \varphi_1(\gamma_2 )
\varphi_2 (\gamma_1) \varphi_2 (\gamma_2)
\\  & = &
\varphi_1(\gamma_1) \varphi_2 (\gamma_1) \varphi_1(\gamma_2 )
\varphi_2 (\gamma_2)
\\  & = &
  (\varphi_1\varphi_2) (\gamma_1)   (\varphi_1\varphi_2) (\gamma_2).
\end{eqnarray*}
The identity $e$ in $\Hom(\pi, G)$ is the trivial representation. Observe finally that this group is abelian since $ \varphi_1 (\gamma) \varphi_2   (\gamma) =  \varphi_2(\gamma) \varphi_1(\gamma)$. \\
\end{proof}

Recall that the  abelianized group of $\pi$ is the abelian group
 $$Ab(\pi) = \pi/[\pi,\pi] .$$
 Another useful remark is that if $G$ is abelian, then
$$\Hom(\pi' \times \pi'', G) = \Hom(\pi', G)\times \Hom(\pi'', G)$$
for any groups $ \pi' ,\pi''$.

 \medskip

Assume now that $\pi$ is a finitely generated group. $Ab(\pi)$
is then an abelian group of finite type, hence
 $$Ab(\pi) = \pi/[\pi,\pi] = \mathbb{Z}^r\oplus F$$
where $F$ is a finite abelian group (the torsion) and $r\in \mathbb{N}$ is the \emph{abelian rank} of $\pi$.

We obviously have $\Hom(\pi, G)= \Hom(Ab(\pi), G)$ and it is clear that all representation varieties of a  finitely generated group $\pi$ in an abelian Lie group $G$ can be deduced from the following special cases:

\smallskip

\begin{enumerate}[1)]
  \item $\Hom(\mathbb{Z}, \R)  = \R $;
  \item $\Hom(\mathbb{Z}, U(1))  = U(1)$;
  \item $\Hom(\mathbb{Z}/m\mathbb{Z}, \R)   = 0$;
  \item $\Hom(\mathbb{Z}/m\mathbb{Z}, U(1))  = \{z \in \C \tq z^m = 1 \}$.
\end{enumerate}

\medskip 

For instance, if $\pi$ is the free group on $s$ generators, then
$Ab(\pi) = \mathbb{Z}^{s}$. Thus
$\Hom(\pi, \R)=\R^{s}$ and $\Hom(\pi, U(1))= \mathbb{T}^{s}$.

Another simple  example, with torsion, is  the group $\pi' = \langle a,b,c \tq
 [a,b]  = c^m \rangle$.  We have
$Ab(\pi') = \mathbb{Z}^2 \oplus \mathbb{Z}/m\mathbb{Z}$, therefore
$\Hom(\pi', \R)=\R^{2}$ and
$$
\Hom(\pi', U(1))= \mathbb{T}^{2} \oplus \{ e^{2ki \pi /m} \tq m=0,1,\dots, m-1\} .
$$

\subsection{Representations in $\SE(2)$}

We denote by $SE(2) = \text{Iso}^+(\mathbb{R}^2)$ the  group of orientation preserving isometries of the  euclidean plane.

\medskip

We may identify the euclidean plane with the complex line $\mathbb{C}$: any $g\in SE(2)$ can then be written as $g(z) = u\cdot z + v$ where $u\in U(1)\subset \mathbb{C}^*$ and $v\in \mathbb{C}$. We thus identify
$\SE(2)$ with the subgroup  of $GL_2(\mathbb{C})$
consisting of matrices of the form
$$
 \SE(2) =\left\{ \left. \left(\begin{array}{cc}u & v \\0 & 1\end{array}\right) \right| u,v \in \mathbb{C}, \ |u|=1\right\}.
$$
In particular  $\SE(2)$ is a semidirect product $U(1)  {\rtimes} \mathbb{C}$ and any representation $\varphi \in \Hom(\pi,SE(2))$ can be written as
\begin{equation}\label{dec.phi}
\varphi  =
\left(\begin{array}{cc} \rho_{\varphi} &\tau_{\varphi} \\0&1\end{array} \right)
\end{equation}
where $\rho_{\varphi} : \pi \to U(1)$ and $\tau_{\varphi} : \pi \to \C$. Observe the following:

\begin{lemma}
The map $\rho_{\varphi} : \pi \to U(1)$  is a group homomorphism. It only depends  on the conjugacy class of $\varphi$.
\end{lemma}
The proof is elementary.

\begin{definition}
The homomorphism  $\rho_{\varphi} : \pi \to U(1)$ is the \emph{character}\index{character} of the
representation class  $\varphi \in \Hom(\pi,SE(2))$.
\end{definition}

\medskip

\textbf{Remark}: In the literature on group representations,
the  character  $\chi_{\varphi} : \pi \to K$  of a representation $\varphi \in \GL_n(K)$ is classically defined to be the trace of the representation. The two notions of characters are equivalent as shown by  the formula
$$\chi_{\varphi} = \Tr \varphi = 1+\rho_{\varphi}.$$

\medskip

Any homomorphism $\rho \in \Hom(\pi,U(1))$ defines a structure of $\pi -$module on $\mathbb{C}$. We will denote by $ \mathbb{C}_{\rho}$ this $\pi -$module, and we have:
\begin{proposition}
Given any pair of maps
$\rho : \pi \to U(1)$ and
 $\tau  : \pi \to \C$,  the map $\varphi : \pi \to SE(2)$ given by (\ref{dec.phi}) is
a group homomorphism if and only if
$\rho \in \Hom(\pi,U(1))$  and $\tau \in Z^1(\pi, \mathbb{C}_{\rho})$.
\end{proposition}

\begin{proof}
 Suppose that  $\varphi : \pi \to SE(2)$  is  given by (\ref{dec.phi}). Then we have
 $$
 \varphi(\gamma_1\gamma_2) =
   \left(\begin{array}{cc} \rho(\gamma_1\gamma_2) &
  \tau (\gamma_1\gamma_2)\\0&1\end{array} \right)
 $$
and
 \begin{eqnarray*}
 \varphi(\gamma_1)\varphi(\gamma_2)
  & = &
    \left(\begin{array}{cc} \rho(\gamma_1)&
  \tau (\gamma_1)\\0&1\end{array} \right)
  \,
      \left(\begin{array}{cc} \rho(\gamma_2)&
  \tau (\gamma_2)\\0&1\end{array} \right)
   \\ & = &
    \left(\begin{array}{cc} \rho(\gamma_1)\rho(\gamma_2) &
  \tau (\gamma_1) + \rho(\gamma_1)\tau(\gamma_2)\\0&1\end{array} \right).
\end{eqnarray*}
It follows that $\varphi$ is  a group homomorphism (i.e. $\varphi(\gamma_1\gamma_2) = \varphi(\gamma_1)\varphi(\gamma_2)$)
if and only if
$$\rho(\gamma_1\gamma_2) = \rho(\gamma_1)\rho(\gamma_2)$$
and
$$
\tau (\gamma_1\gamma_2) =
\tau (\gamma_1) + \rho(\gamma_1)\tau(\gamma_2).
$$
In other words $\varphi$ is  a   homomorphism if and only if
$\rho : \pi \to U(1)$ is  a   homomorphism and $\tau$ is a 1-cocycle
in the corresponding $\pi$-module $ \mathbb{C}_{\rho}$.\\
\end{proof}

This Proposition says that
the map from  $\Hom(\pi  , \SE(2))$ to the set
$$
 \{ (\rho , \tau) \tq  \rho \in \Hom(\pi, U(1)) \ \mathrm{and} \  \tau\in
 Z^1(\pi, \mathbb{C}_{\rho}) \}
$$
given by
$
 \varphi \to (\rho_{\varphi} , \tau_{\varphi}),
$
is a bijection. In particular we have

\begin{corollary}
 If $\pi$ is the free group on $s$ generators, then
$$ \Hom(\pi  , \SE(2)) \simeq \mathbb{T}^{s} \times \mathbb{C}^s.$$
\end{corollary}

\begin{proof}
This follows from equation (\ref{eq.Z1}) and the fact that
$ \Hom(\pi  , U(1)) = \mathbb{T}^{s}$.

\end{proof}

\subsection{Conjugation by similarities}

Recall that a \emph{similarity} in the plane is
the composition of an isometry with a homothety.

\medskip

Once we identify the euclidean plane with the complex line
$\mathbb{C}$, any similarity $g\in \Sim(2)$ can be writen as $g(z) =
a\cdot z + b$ where $a\in \mathbb{C}^*$ and $b\in \mathbb{C}$. We
thus identify $\Sim(2)$ with the following  subgroup  of
$GL_2(\mathbb{C})$:
$$
 \Sim(2) =\left\{ \left. \left(\begin{array}{cc}a & b \\0 & 1\end{array}\right) \right| a,b \in \mathbb{C}, \ a\neq 0 \right\}.
$$
In particular  we have
$$
 \Sim(2) = \mathbb{R}_+\rtimes \  \SE(2)= \mathbb{C}^* {\rtimes} \ \mathbb{C}.
$$

\begin{definition}
Two representations $\varphi_1, \varphi_2: \pi \to \SE(2)$ are
\emph{similar}\index{similar representations} if they are
conjugate modulo a similarity.
\end{definition}

\begin{proposition}\label{prop.repcohom}
Given a homomorphism $\rho : \pi \to U(1)$ and two cocycles
$\tau_1,\tau_2  \in Z^1(\pi , \C_{\rho})$, then the representations
\begin{equation}\label{}
\varphi_1  =
\left(\begin{array}{cc} \rho &\tau_1 \\0&1\end{array} \right)
\quad and \quad
\varphi_2  =
\left(\begin{array}{cc} \rho &\tau_2 \\0&1\end{array} \right)
\end{equation}
are similar  if and only if there exists a complex number $a\in
\mathbb{C}^*$  such that
$$\tau_2 = a \tau_1 \in H^1(\pi, \mathbb{C}_{\rho}).$$
\end{proposition}

\medskip

\begin{proof}
The homomorphisms $\varphi_1$ and $\varphi_2$ are similar if and
only if there exists
$$
 g = \left(\begin{array}{cc} a & b \\0&1\end{array} \right)
 \in \Sim(2)
$$
such that $\varphi_2 = g\varphi_1 g^{-1}$, i.e.
 \begin{eqnarray*}
\left(\begin{array}{cc} \rho &\tau_2 \\0&1\end{array} \right) & =
& \left(\begin{array}{cc} a & b \\0&1\end{array} \right)
\left(\begin{array}{cc} \rho &\tau_1 \\0&1\end{array} \right)
\left(\begin{array}{cc} 1/a & -b/a \\0&1\end{array} \right)
\\ & = &
 \left(\begin{array}{cc} \rho & a\tau_1 + b - \rho b \\0&1\end{array} \right).
\end{eqnarray*}
This shows that
$$ \tau_2 - a\tau_1 = b\, (1 - \rho ) \in  B^1(\pi, \mathbb{C}_{\rho}).$$
\end{proof}

\bigskip

For any homomorphism $\varphi : \pi \to \SE(2)$ and any $\lambda\in \R_+$,
we can define a new  homomorphism $\lambda\cdot\varphi : \pi \to \SE(2)$
by

$$
\lambda\cdot \varphi  =
\left(\begin{array}{cc} \rho &\lambda \tau \\0&1\end{array} \right),
$$
This formula defines an action of the  multiplicative group $\R_+$  on
$\mathcal{R }(\pi  , \SE(2))$, and we denote the quotient by
$$
\mathcal{SR}(\pi  , \SE(2)) = \mathcal{R }(\pi  , \SE(2))/\R_+ .
$$
It follows directly from the definition that
$$
\mathcal{SR}(\pi  , \SE(2))  =
\Hom (\pi , \SE(2)) /\Sim(2).
$$
where $\Sim(2)$ acts by conjugation on
$\Hom (\pi , \SE(2))$.

Let us also define
$$
 \mathcal{SR}^{\mbox{reg}} =
  \{ [\varphi] = [\rho, \tau]  \in \mathcal{SR }(\pi  , \SE(2))
 \tq  \rho_{\varphi} \neq id \ \mathrm{and }\ \tau \neq 0
 \}
 $$

\begin{corollary} \label{cor.srreg}
If  $\pi$ is a free group on $s$ generators, then
$$
 \mathcal{SR}^{reg}  \simeq
(\mathbb{T}^{s} \setminus \{ id \}) \times \mathbb{CP}^{s -2}
$$
\end{corollary}

\begin{proof}
This is an immediate consequence of the previous  results, in particular Proposition \ref{compute.H1} and \ref{prop.repcohom}.

\end{proof}


\section{Deformation Theory}
\label{sec.teich}


\subsection{The Moduli and Teichm\"uller spaces}

The  moduli space of $\Sigma_{g,n}$ is the quotient of the space
of conformal structures on $\Sigma_{g,n}$ by the pure
diffeomorphism group. Let us be more specific: recall first that
a conformal structure is an equivalence class of  smooth Riemannian
metric $m$ on $\Sigma = \Sigma_{g,n}$, where  two Riemannian metrics
$m_1,m_2$ are equivalent if and only if there exists a  function
$u : \Sigma \to \R$ such that
$$
 m_2 = e^{2u} m_1.
$$
We denote by $\Met(\Sigma)$ the space of all smooth Riemanian metrics on $\Sigma$ endowed with
its natural $C^{\infty}$ topology and by
$$
 \Conf(\Sigma) =  \Met(\Sigma)/ C^{\infty}(\Sigma)
$$
the space of conformal structures. 
We then define the \emph{moduli space}\index{moduli space} of $\Sigma_{g,n}$ to be the quotient
$$
 \mathcal{M}_{g,n} =  \Conf(\Sigma)/\PDiff^+_{g,n}.
$$

\medskip

A point $\mu \in \mathcal{M}_{g,n}$ is concretely represented by
a Riemannian metric $m$ on $\Sigma$, and  two Riemannian
metrics $m_1,m_2$  represent the same modulus  point $\mu$ if and only if
there exists a smooth function $u : \Sigma \to \R$ and a
diffeomorphism $h \in \PDiff^+_{g,n}$ such that $m_2 = e^{2u}
h^*(m_1)$.

\medskip

{\footnotesize
A remark about the smoothness: By definition a point $\mu$ in the moduli space is represented by a smooth metric.
In particular,  the puntures play no role in the definition of the spaces  $\Met(\Sigma)$ and $\Conf(\Sigma)$ (but
they do in  the definition of the moduli space  $\mathcal{M}_{g,n}$).
However, since only the conformal class of the metric matters, one may also represent $\mu$ by a singular metric $m$ as long as this metric is conformally equivalent to a smooth one. In particular we can (and will) represent a point in $\mathcal{M}_{g,n}$ by a metric having conical singularities at the punctures of $\Sigma_{g,n}$, see Remark \ref{sing_supprimable}
}%

\medskip

The moduli space is a complicated object, and it is useful to
also introduce the simpler space obtaind by considering isotopy
classes of  conformal structures on $\Sigma_{g,n}$ instead of
isomorphism classes: this is the   \emph{Teichm\"uller
space}\index{Teichm\"uller space} defined as
$$
 \mathcal{T}_{g,n}  =  \Conf(\Sigma_{g,n})/\Diff^0_{g,n}.
$$

\bigskip

Let us list some of the basic  facts about these spaces:

\begin{enumerate}[(1)]
  \item The Teichm\"uller space $\mathcal{T}_{g,n}$ is a real analytic variety in a natural way. If $3g-3 +n >0$, then it is isomorphic to $\R^{6g - 6 + 2n}$. This space has also  a natural complex  structure.
  \item The pure mapping class group  $\PMod^+_{g,n}$ acts properly and discontinuously on $\mathcal{T}_{g,n}$.
  \item The moduli space is the quotient $\mathcal{M}_{g,n} = \mathcal{T}_{g,n}/\PMod_{g,n}^+$. It is thus a good orbifold  of dimension $6g - 6 + 2n$ with fundamental group $\pi_1(\mathcal{M}_{g,n}) = \PMod_{g,n}$.
\item There exists a torsion free subgroup $M_0 \subset \PMod^+_{g,n}$ of finite index acting freely on $\mathcal{T}_{g,n}$. The quotient map
$\mathcal{T}_{g,n}/M_0$ is a non singular analytic manifold which
is a finite cover of the orbifold $\mathcal{M}_{g,n}$.
\end{enumerate}

\begin{proof}
 Statement (1) is explained in any textbook on Teichm\"uller theory such as \cite{abikoff}. Statement  (2) was first proved by S. Kravetz \cite{kravetz}, see also \cite{abikoff}. (3) is a consequence  of (1) and (2) and the last statement is discussed in \cite[\S 5.4]{ivanov02}.

\end{proof}

\subsection{The deformation space of piecewise flat metrics}

Let us  denote by $\mathcal{E}_{g,n}$ the  set of all  flat
metrics on $\Sigma_{g,n}$  with possible conical singularities at
the punctures (it is not empty since we have assumed  
$2g+n-2> 0$). To any flat metric  $m\in  \mathcal{E}_{g,n}$, we
associate the following basic invariants~:  Its conformal class  
$[m] \in \Conf({\Sigma})$, its area $A = A(m) > 0$ and the order $\b_i
> -1$ of $m$ at the point $p_i$.
\begin{theorem}
The map
\begin{align*}
   \mathcal{E}_{g,n}   &\to \Conf (\Sigma) \times \R_+^n \\
    m  &  \mapsto ([m], (s_1,...,s_n) ),
\end{align*}
where $s_i = A\, (1+\b_i) > 0$, is a bijection.
\end{theorem}

\begin{proof}
 This is just a reformulation of Theorem  \ref{th.srsc}.

\end{proof}

\medskip

\begin{definition}
 We will endow the set $\mathcal{E}_{g,n}$ with the topology for which this map
 is a homeomorphism.
\end{definition}

\medskip

A metric $m_2\in  \mathcal{E}_{g,n}$ is said to be a
\emph{deformation} of the metric $m_1\in  \mathcal{E}_{g,n}$ if the two  metrics
 differ by a homothety and an isotopy fixing the punctures,
i.e. if there exists   $h\in \PDiff^0_{g,n}$ and $\lambda > 0$ such that
$m_2 = \lambda \, h^* (m_1)$. We denote by $\mathcal{DE}_{g,n}$
the deformation space of flat metrics on $\Sigma_{g,n}$  with
possible conical singularities at the punctures:
$$
 \mathcal{DE}_{g,n} = \mathcal{E}_{g,n}/(\R_+\times \PDiff^0_{g,n}).
$$
\begin{corollary}
This space is homeomorphic to $\R^{6g+3n-7}$. In fact we have the
following canonical identification:
$$
  \mathcal{DE}_{g,n} = \mathcal{T}_{g,n}\times \Delta,
$$
where  $\mathcal{T}_{g,n}$ is the Teichm\"uller space and
$\Delta \subset \R^n$ is  defined by
$$
\Delta =  \{ \vec{\b} = (\b_1,...\b_n ) \in \R^n \tq  \b_i > -1
\ and \ \sum_i \b_i = 2g-2\}.
$$
\end{corollary}

\medskip

Let us fix an element $\vec{\b}  = (\b_1,...\b_n)\in \Delta$ and 
denote by $\mathcal{E}_{g,n}(\vec{\b})$ the space of singular
flat metrics with a conical singularity of order $\b_i$ at $p_i$
($i=1,...,n)$. We also introduce the corresponding deformation
space: $\mathcal{DE}_{g,n}(\vec{\b}) =
\mathcal{E}_{g,n}(\vec{\b})/(\R_+\times \PDiff_{g,n})$. The
previous corollary gives us the identification
$$
 \mathcal{DE}_{g,n}(\vec{\b})= \mathcal{T}_{g,n}.
$$

\subsection{Revisiting the development and the holonomy}

Consider the punctured surface $\Sigma_{g,n} =
\widehat{\mathbb{U}}/\Gamma$ as in section
\ref{subs.uniformization}, and fix a flat metric  $m_0$ with  conical singularity of order $\b_j$ at $p_j$  ($j =1,\dots , n$).  If  $f_0$ is a germ of an isometry  near a point 
$\tilde{z}_0$,  to the euclidean plane (identified with
$\mathbb{C}$), then we obtain a map
 $f:\mathbb{U} \to \mathbb{C}$ by analytic continuation from $f_0$.
This map is a local isometry for the metric $m_0$ on $\mathbb{U}$
and the canonical metric on  $\mathbb{C}$ (indeed, the set of points 
where a map $f$ between two flat surfaces is an isometry is easily seen to be both open and closed). The map $f$  extends by
continuity to $\widehat{\mathbb{U}}$.  The resulting map $f :
\widehat{\mathbb{U}} \to  \mathbb{C}$ is the developing map which
we already met in Section \ref{sec.dev}. The associated holonomy
is the unique  homomorphism $\varphi : \Gamma \to \SE(2)$ such
that
$$
  f(\gamma u) = \varphi (\gamma) f (u).
$$

\begin{theorem} The following properties of the development and its
holonomy are satisfied:
\begin{enumerate}
  \item $f :  \widehat{\mathbb{U}} \to  \mathbb{C}$ is surjective;
  \item $f(\gamma u) = \varphi (\gamma) f (u)$ for all $\gamma \in \Gamma$;
  \item for any $y\in P^{-1}(p_i) \subset \Upsilon$, the isometry
  $\varphi (c_y)\in \SE(2)$ is a rotation of angle $\theta_i = 2\pi(\b_i+1)$;
  \item if $\b_i$ is not an integer, then $f(y)$ is the unique fixed point of
  $\varphi (c_y)$.
\end{enumerate}
\end{theorem}

\begin{proof}
The first assertion has been proved in Proposition
\ref{dev.surjective}, the second is the definition of the holonomy
and the last two assertions are contained in Proposition
\ref{local.dev}.

\end{proof}

\begin{corollary}\label{cor.holmapdev}
If $\b_i \not\in \mathbb{Z}$ for any $i=1,...,n$, then the
restriction of $f$ to the set $\Upsilon$ is determined by the
holonomy.
\end{corollary}

\begin{proof}
 Fix $y\in \Upsilon$ and let $c_y \in \Gamma$ be the corresponding group element given by the map (\ref{Ydevmap}).
Then $f(y)\in \mathbb{C}$ is the fixed point of the rotation
$\varphi(c_y)\in \SE(2)$. This fixed point is given explicitly by 
\begin{equation}\label{imageY}
 f(y)= \frac{\tau_y}{(1-\rho_y)},
\end{equation}
where $\rho_y \in U(1)$ is the rotation part  and 
$\tau_y\in \mathbb{C}$ is the translation part of $\varphi(c_y)$.
\end{proof}
 
\bigskip

\begin{theorem}\label{mainth.a}
\emph{(A) }There is a well defined   map
$$
 hol : \mathcal{DE}_{g,n} \to \mathcal{SR}(\pi_{g,n},\SE(2)),
$$
such that $hol([m])$ is the conjugacy class of the  holonomy homomorphism $\varphi _m : \pi_{g,n}\to \SE(2)$.

\medskip

\emph{(B)} The map $hol : \mathcal{DE}_{g,n} \to \mathcal{SR}(\pi_{g,n},\SE(2))$ is continuous.

\medskip

\emph{(C)} There are natural actions of $\PMod_{g,n}$ on
$\mathcal{DE}_{g,n}$ and  $\POut(\pi_{g,n})$ on $\mathcal{SR}(\pi  , \SE(2))$, and the map $hol$ is
$\Phi$-equivariant where $\Phi$  is the Dehn-Nielsen-Baer isomorphism (\ref{eq.DNB}).

\medskip

\emph{(D)} The map   $hol$ is locally injective.
\end{theorem}

\textbf{Remarks.}  1.) The map $hol : \mathcal{DE}_{g,n} \to \mathcal{SR}(\pi_{g,n},\SE(2))$ is called the \emph{holonomy mapping}. \\
2.) A more elaborate investigation would show that the holonomy mapping is
in fact real analytic, see \cite{veech93}. The proof below is perhaps not
optimal from the point of view of rigour, but we have tried to  emphasize the geometric point of view.

\begin{proof}
(A) To any flat metric $m$ on $\Sigma_{g,n}$ with conical singularities at the
punctures, we have associated a holonomy homomorphism $\varphi _m : \pi_{g,n}\to \SE(2)$
which depends on the choice of a developing map $f_m$, but changing the  developing map does not affect the  conjugacy class of $\varphi _m$ (see Proposition \ref{prop.conjhol}).
On the other hand, it is clear that if two flat metrics $m,m'$ on $\Sigma_{g,n}$ are similar, the associated holonomies  $\varphi _m, \varphi _{m'}$ are also similar. In short, to any deformation class of  flat metric $[m] \in \mathcal{DE}_{g,n}$ with conical singularities on $\Sigma_{g,n}$ we associate 
a well defined element $[\varphi_m] = hol(m) \in \mathcal{SR}(\pi_{g,n},\SE(2))$.

\medskip

(B) The developing map $f_m$ is not uniquely associated to a flat metric $m$, but it is well defined modulo $SE(2)$ (two developing maps for the same metric differ by postcomposition with an isometry). 
The $SE(2)$ orbit of the developing map   $f_m$ varies continuously with the metric $m$
and therefore it is also the case for the associated holonomy class. Hence the map $hol : \mathcal{DE}_{g,n} \to \mathcal{SR}(\pi_{g,n},\SE(2))$ is continuous.
 
 \medskip

(C) Any diffeomorphism $h$ of $\Sigma_{g,n}$ fixing the punctures acts on 
$\mathcal{E}_{g,n}$ by pulling back the metric ($m \mapsto h^*m$). If $h$ is isotopic
to the identity, it acts trivially on $\mathcal{DE}_{g,n}$; we thus have a well defined action
of $\PMod_{g,n}$ on $\mathcal{DE}_{g,n}$.

Similarly, any automorphism of $\pi_{g,n}$ acts on $\Hom(\pi_{g,n},\SE(2))$, and inner automorphisms act trivially on  the representation spaces $\mathcal{R}(\pi_{g,n} , \SE(2))$
and $\mathcal{SR}(\pi_{g,n}  , \SE(2))$. We thus have a natural action of  $\POut(\pi_{g,n})$
on these spaces. It is clear from the construction of the isomorphism
$\Phi : \PMod_{g,n} \stackrel{\sim}{\longrightarrow} \POut(\pi_{g,n})$
(see the proof of Lemma \ref{lem.DNB}) that the map $hol$ is equivariant.

\medskip

(D) To prove the local injectivity of $hol$, we consider two nearby flat metrics $m,m'$ with conical singularities on  $\Sigma_{g,n} = \widehat{\mathbb{U}}/\Gamma$
and we assume that they have the same holonomy $\varphi$. Since the holonomy around a conical singularity $p_i$ is a rotation of angle $\theta_i = 2\pi(\b_i+1)$, it is clear that both metrics $m$ and $m'$ have the same singularity order (the holonomy only controls the cone angle modulo $2\pi$, but since $m$ and $m'$ are nearby metrics, 
they actually have equal cone angles).

It follows that both metrics are isometric near the singularities: we can thus find an isotopy $h_1$ of the surface such that $m=h_1^*m'$ near the singularities. Hence we can simply assume without loss of generality  that  $m=m'$ near the singularities; it is therefore possible to divide the surface in $n+1$ parts
$$\Sigma_{g,n} = D\cup E_1 \cup \cdots \cup E_n ,$$
where $D\subset \Sigma'$ is a compact region and $E_i$ is a neighbourhood of the puncture $p_i$ such that $m=m'$ on $E_i$. We also assume that the $E_i$ are pairwise disjoint disks. We denote by $\widehat{E}_i = P^{-1}(E_i) \subset \widehat{\mathbb{U}}$ and $\widehat{D} =P^{-1}(D) \subset \widehat{\mathbb{U}}$ the lifts of $E_i$ and $D$ on the universal branched cover $P : \widehat{\mathbb{U}} \to \Sigma_{g,n}$. The set $\widehat{E} = \cup_i \widehat{E}_i  \subset \widehat{\mathbb{U}}$ 
is a neighbourhood of $\Upsilon = P^{-1}(\{\mathrm{punctures}\})$.

Let $f_m$ and $f_{m'}$ be the developing maps of $m$ and $m'$. By Corollary \ref{cor.holmapdev}, the two maps coincide on $\Upsilon$. Because the two metrics coincide on ${E}_i$,  the map $f_m$ and $f_{m'}$ coincide up to
a rotation on each component of $\widehat{E}_i$; we can thus find a second isotopy $h_2$ of $\Sigma$, which is a rotation near the punctures and  is the identity on $D$ and such that $f_{m'} \circ \widehat{h}_2 = f_{m}$ on $\widehat{E}$.

Replacing $m'$ with $h_2^*m'$, we can thus assume that both  developing maps 
coincide on  $\widehat{E}$.

To any point $x\in \mathbb{U}$, we associate the set 
$$\Lambda (m,m',x) = f_{m}^{-1}(f_{m'}(x)) \subset \widehat{\mathbb{U}}.$$ 
Since $f_m$ is a local diffeomorphism, $\Lambda (m,m',x)$ is
a discrete set. It varies continuously with $m,m'$.

\smallskip

Claim: \emph{If $m$ is close enough to $m'$, then for any point $x\in U$, there exists a unique point $y = Q(x) \in \Lambda (m,m',x) $ which is the nearest point for the
hyperbolic distance. The map $x\mapsto Q(x)$ is $\Gamma$-equivariant.}

\smallskip

Indeed, if  $x\in \widehat{E}$, then the claim is clear: since $f_{m}(x) = f_{m'}(x)$, we
have $Q(x)=x$. For any point, the claim is clear if $m=m'$ (and in this case $Q(x) = x$). 
For points in $\widehat{D}$, and $m'$ close to $m$ the claim follows from the compactness of $D$
and the discreteness and continuity  of $\Lambda (x,m,m')$.

\smallskip

For $t\in [0,1]$, we denote by $Q_t(x)$ the point on the hyperbolic segment $[x,Q(x)]$ such that $d_H(x,Q_t(x)) = t d_H(x,Q(x))$ 
(observe that if  $x\in \widehat{E}$, then $Q_t(x)=x$ for any $t$).  This is a $\Gamma$-equivariant isotopy of $U$ from the identity to $Q$. It extends as the
identity on $\Upsilon$.

We now define an isotopy $h_t : \Sigma \to \Sigma$ by $h_t(x) = P(Q_t (P^{-1}(x)))$. It is a well defined isotopy such that $h_1^*m' = m$, since we clearly have  $f_m = f_{m'}\circ Q$. 

We thus have proved that two metrics with the same holonomy are isotopic
provided they are close enough. In other words, the map $hol$ is locally injective.

 \end{proof}


\section{The Main Theorem}

We are now in position to prove the main result. First recall the
statement:

\begin{theorem} Given a punctured surface $\Sigma_{g,n}$
such that $2g+n-2> 0$ and $\vec{\b}\in \Delta$ such that no $\b_i$
is an integer,  there is a well defined group homomorphism
$$
 \Phi : \PMod_{g,n} \to \mathcal{G} = \Aut (\mathbb{T}^{2g})\times \PGL_{2g+n-2}\mathbb{C},
$$
and a  $\Phi$-equivariant local homeomorphism
$$
\mathcal{H} : \mathcal{T}_{g,n} \to  \Xi = \mathbb{T}^{2g} \times
 \mathbb{CP}^{2g+n-3}.
$$
\end{theorem}

The theorem says that $ \mathcal{M}_{g,n}  =  \mathcal{T}_{g,n} /\PMod_{g,n}$ is a good orbifold with a
$(\mathcal{G}, \Xi)$-structure.

\medskip

\textbf{Proof}  The group homomorphism $\Phi$ is given by the Dehn-Nielsen-Baer
isomorphism and the map $\mathcal{H}$ is essentially given by the holonomy
mapping of the previous theorem. We  divide the proof of the theorem in 5 steps:

\medskip

Recall that holonomy splits in a rotation
part $\rho _m : \pi \to U(1)$ (the character) and a translation
part $\tau _m$. The character depends only on the conjugacy class
of $\varphi _m$.

\medskip

\textbf{Step 1:}  \emph{ $\tau_m$ is not identically zero.}

\medskip

Indeed, if $\tau_m \equiv 0$, then the holonomy group $\varphi_m(\pi_{g,n})$
is a pure rotation group in the plane. This is impossible by Proposition \ref{noboundedorbit}.

\medskip

\textbf{Step 2:} \emph{There is a canonical isomorphism
$$
\left.\begin{array}{ccc} \Hom(\pi_{g,n} ,U(1))  & \simeq  &
\Hom(\pi_{g,0} ,U(1)) \times \Hom(\pi_{0,n},U(1)).  \\
\rho_m & \mapsto & (\rho ' , \rho '')
\end{array}\right.
$$
Furthermore $\rho '' \in  \Hom(\pi_{0,n},U(1))$ is given by
$$\rho ''(c_i) = e^{\theta_i},$$
where $c_i$ is the homotopy class of a loop traveling once around
the puncture $p_i$ and $\theta_i = 2\pi(\b_i +1)$ is the total
angle at the cone point $p_i$.
}

\medskip

This splitting easily follows from the  identity  (\ref{dec.Psurface})  and the fact
that $U(1)$ is abelian.

\medskip

Let us now fix an element $\vec{\b}  = (\b_1,...\b_n)\in \Delta$ and
set
$$
\mathcal{SR}_{\vec{\b}}(\pi_{g,n} , \SE(2)) = \left\{ \varphi \in
\mathcal{SR}(\pi_{g,n} , \SE(2)) \tq \rho ''(c_i) =
e^{\theta_i}, i = 1,...,n \right\}
$$
and
$$
 \mathcal{SR}^{reg}_{\vec{\b}} =
  \mathcal{SR}^{reg} \cap  \mathcal{SR}_{\vec{\b}}.
$$

\textbf{Step 3:} \emph{If at least one $\b_i$ is not an integer, then we have}
$$
\mathcal{SR}^{reg}_{\vec{\b}}(\pi_{g,n} , \SE(2)) \simeq \Xi =
\mathbb{T}^{2g} \times \mathbb{CP}^{2g+n -3}
$$

\medskip

Indeed, it follows from Step 2  and the results of Section \ref{sec.rep} that 
 any $\varphi  \in \mathcal{SR}^{reg}_{\vec{\b}}$
is characterized by $\rho' \in \Hom(\pi_{g,0} ,U(1)) \simeq \mathbb{T}^{2g}$ and
the projective class of $\tau \in H^1(\pi_{g,n},\mathbb{C}_{\rho}) \simeq \mathbb{C}^{2g+n-2}$ \
(because  $\pi_{g,n}$ is isomorphic to the free group on $s=2g+n-1$
generators).

\medskip

\textbf{Step 4:} \emph{The group $\POut(\pi_{g,n})$ acts naturally on 
$\Xi$ and thus 
we have a natural homomorphism 
$ \Phi : \PMod_{g,n} \to \mathcal{G} = \Aut (\mathbb{T}^{2g})\times \PGL_{2g+n-2}\mathbb{C}$.}

\medskip
 
This is clear from Step 3 and part (C) of Theorem \ref{mainth.a}.

\medskip

\textbf{Step 5:} \emph{The map $\mathcal{H}$ given by  the composition:
$$
\mathcal{T}_{g,n} \stackrel{\sim}{\to}
\mathcal{DE}_{g,n}(\vec{\b}) \stackrel{hol}{\to}
\mathcal{SR}^{reg}_{\vec{\b}}\stackrel{\sim}{\to}\, \Xi
$$
is well defined,  continuous, locally injective and $\Phi$-equivariant.
}

\medskip

Indeed, the fact that no $\b_i$ is integer, together with Step 1, implies that  $hol : \mathcal{DE}_{g,n} \to
\mathcal{SR}(\pi , \SE(2))$ maps $\mathcal{DE}_{g,n}(\vec{\b})$
into $\mathcal{SR}^{reg}_{\vec{\b}}$. \ 
The map $\mathcal{H} : \mathcal{T}_{g,n} \to
\Xi$ is therefore well defined.
It follows from Theorem \ref{mainth.a} that $\mathcal{H}$ is
continuous, locally injective and $\Phi$-equivariant.

It remains only to show that $\mathcal{H}$ is a local homeomorphism, but since $\mathcal{T}_{g,n}$ and
$\Xi$ are both manifolds of dimension $6g-6+ 2n$, the conclusion
follows from  Brouwer's Theorem on invariance of dimension.

\qed


\subsection{The case of the sphere}

Suppose that  $g = 0$, i.e. $\Sigma$ is a sphere,  choose $n$
numbers ($n \ge 2$) $\b_1,\b_2,\dots ,\b_n$ such that $2 + \sum_i
\b_i = 0$, and denote by $\mathcal{M }$ the space of flat metrics
on $S^2$ having $n$ conical singularities of order
$\b_1,\b_2,\dots ,\b_n$.

\medskip

Such a metric $m \in \mathcal{M}$ can be uniformized as follows :
identify $\Sigma$ with $\mathbb{C} \cup \infty$, and write $m$ as
$$m = C\cdot \prod_{i=1}^n |z-p_i|^{2\b_i}|dz|^2 \, ,$$
where $p_1,p_2,\dots , p_n$ is the set of conical singularities
and
 $C$ is a positive constant representing a dilation factor.
 It is easy to see that $\mathcal{M}$ is homeomorphic to the quotient
$$\{ (p_1,p_2,\dots , p_n) \in (\mathbb{C}  \cup \infty )^n :
p_i \ne p_j \text{ if } i \ne j \}/PSL_2(\mathbb{C}) \, ,$$
$\mathcal{M}$ is thus a complex manifold of dimension $n-3$, its
fundamental group is the  pure braid group $PB_n$.

\medskip

Applying the main theorem,  we obtain a  representation
$$
 \Phi : \PMod_{0;n}  = PB_n \to  \PGL_{n-2}(\mathbb{C})
$$
and a $\Phi$ equivariant, local homeomorphism
$$
\mathcal{H} : \mathcal{T}_{0,n} \to
 \mathbb{CP}^{n-3}.
$$

\medskip

In fact, a finer analysis shows that the image of $\Phi$ is
contained in $PU(1,n-3) \subset  \PGL_{n-2}(\mathbb{C})$.
Furthermore, when the orders satisfy some  arithmetical
conditions, the image of $\Phi$ is a lattice in $PU(1,n-3)$:

\begin{theorem}
 Assume that $-1<\b_1,\b_2, ...,\b_n < 0$, $\sum_i \b_i = -2$ and suppose that
 \begin{equation}\label{cond.DM}
\b_i + \b_j > -1\,  \Rightarrow \, (1+\b_j+\b_i)^{-1} \in
\mathbb{N},
\end{equation}
then $\Phi(PB_n)$ is a lattice in $PU(1,n-3)$. \\
These lattices are quotients of the braid group. Some of them are
non arithmetic.
\end{theorem}

This Theorem was first proved by Schwartz (1873) for $n=4$ and by
Picard (1888) for $n=5$ in their study of the monodromy of the
hypergeometric equations. It has been generalized for any $n$ by
P. Deligne and G. Mostow in 1986, see  \cite{DM}. These authors
use the cohomology with coefficients in flat vector bundle on an
algebraic curve.

In the paper  \cite{thurston98}, W. Thurston obtain the same
result by studying a deformation space of piecewise flat
triangulations on the sphere (this nice paper is  a 1987 preprint
of W. Thurston, which has been rewritten and appeared in
electronic form in 1998). It is worthwile to quote also the related
papers \cite{BG,fillastre,mostow,sadayoshi}.

\medskip 

Our approach can be seen as a bridge between the approach of
Thurston and that of Deligne-Mostow.

\medskip

Observe that the moduli space $\mathcal{M} =
\mathcal{T}_{0,n}/PB_n$ carries a complex hyperbolic metric
(depending upon the choice of the $\b_i$'s). It is not complete
as a Riemannian manifold and it carries  a natural completion
$\overline{\mathcal{M}}$. Thurston shows that
$\overline{\mathcal{M}}$ is a complex hyperbolic manifold with
singularities of conical type. This cone-manifold has finite volume.

Furthermore, when the $\b_i$'s satisfy the condition (\ref{cond.DM}),
then $\overline{\mathcal{M}}$  is an orbifold. It is thus
possible to  construct  complete complex hyperbolic orbifolds
$\overline{\mathcal{M}}$ of finite volume.

 

\printindex

\end{document}